\documentclass[twocolumn]{autart}
\usepackage{}    
\usepackage{subfigure}
\usepackage{amsfonts}
\usepackage{mathrsfs}
\usepackage{amsmath,amssymb}
\usepackage{graphicx}          
\def\dref#1{(\ref{#1})}

\begin{document}

\begin{frontmatter}

\title{Distributed Robust Consensus Control of Multi-agent Systems
with Heterogeneous Matching Uncertainties} 


\author[li]{Zhongkui Li\corauthref{cor}}\ead{zhongkli@gmail.com},    
\author[li]{Zhisheng Duan}\ead{duanzs@pku.edu.cn},  
\author[lewis]{Frank L. Lewis}\ead{lewis@uta.edu}  
\corauth[cor]{Corresponding author. }

\address[li]{State Key Laboratory for Turbulence
and Complex Systems,
Department of Mechanics and Aerospace Engineering, College of Engineering, Peking University, Beijing 100871, China} %
\address[lewis]{The Automation and Robotics Research Institute, The University of Texas at Arlington, Fort Worth, TX 76118-7115, USA}


\begin{abstract}                          
 This paper considers the
distributed consensus problem of
linear multi-agent systems subject to different matching uncertainties
for both the cases
without and with a leader of bounded unknown control input.
Due to the existence of nonidentical uncertainties,
the multi-agent systems discussed in this paper are essentially
heterogeneous. For the case where the communication graph is undirected and connected,
a distributed continuous
static consensus protocol based on the relative state
information is first designed,
under which the consensus error is uniformly
ultimately bounded and exponentially converges
to a small adjustable residual set.
A fully distributed adaptive consensus protocol is then designed,
which, contrary to the static protocol, relies on
neither the eigenvalues of the Laplacian matrix
nor the upper bounds of the uncertainties.
For the case where there exists a leader whose control
input is unknown and bounded, distributed static and adaptive consensus
protocols are proposed to ensure
the boundedness of the consensus error.
It is also shown that the proposed protocols can be redesigned
so as to ensure the boundedness of the consensus error in the presence of
bounded external disturbances which do not satisfy the
matching condition.
A sufficient condition for the existence of the proposed
protocols is that each agent is stabilizable.
\end{abstract}

\end{frontmatter}

\section{Introduction}

Cooperative control of a network of autonomous agents
has been an emerging research direction and
attracted a lot of attention from many scientific communities,
especially the systems and control community.
A group of
autonomous agents, by coordinating with each other via communication
or sensing networks, can perform certain challenging tasks
which cannot be well accomplished by a single agent.
Cooperative control of multi-agent systems has potential
applications in broad areas including spacecraft formation flying, sensor networks,
and cooperative surveillance \cite{olfati-saber2007consensus,ren2007information}.
In the area of cooperative control,
consensus is an important and fundamental problem,
which means to develop
distributed control policies using only local information
to ensure that the agents reach an agreement on certain quantities of interest.

Two pioneering works on consensus are \cite{jadbabaie2003coordination} and \cite{olfati-saber2004consensus}.
A theoretical explanation
is provided in \cite{jadbabaie2003coordination} for the alignment behavior observed in the Vicsek
model \cite{vicsek1995novel} and a general framework of the consensus
problem for networks of integrators is proposed in \cite{olfati-saber2004consensus}.
Since then, the consensus problem has been extensively
studied by various scholars from different perspectives;
see \cite{olfati-saber2007consensus,ren2007information,ren2005consensus,hong2008distributed,li2010distributed,li2010consensus,ni2010leader,zhang2011optimal,you2011network}
and references therein. Existing
consensus algorithms can be roughly categorized into two classes,
namely, consensus without a leader
(i.e., leaderless consensus) and
consensus with a leader.
The latter is also called leader-follower
consensus or distributed tracking.
In \cite{ren2005consensus}, a sufficient
condition is derived to achieve consensus for multi-agent systems
with jointly connected communication graphs.
The authors in
\cite{hong2008distributed} design a distributed
neighbor-based estimator to track an active leader. Distributed
tracking algorithms are proposed in \cite{ren2010tracking} and
\cite{cao2009distributed} for a network of agents with first-order
dynamics.
Consensus of networks of double- and high-order integrators is studied
in \cite{ren2007distributed,jiang2009consensus}.
Consensus algorithms are designed in \cite{li2010distributed,carli2009quantized} for multi-agent
systems with quantized communication links.
The authors in
\cite{mei2011tracking} address a distributed tracking
problem for multiple Euler-Lagrange systems with a dynamic leader.
The consensus problem of multi-agent
systems with general discrete- and continuous-time
linear dynamics is studied
in \cite{li2010consensus,ni2010leader,zhang2011optimal,you2011network,tuna2009conditions,seo2009consensus,ma2010necessary}.
It is worth noting that the design of the consensus protocols
in \cite{li2010consensus,ni2010leader,zhang2011optimal,seo2009consensus,ma2010necessary}
requires the knowledge of the eigenvalues of the Laplacian matrix of the communication graph,
which is actually global information. To overcome this limitation,
distributed adaptive consensus protocols
are proposed in \cite{li2011adaptive,li2012adaptiveauto}.
For the case where there exists a leader with possibly
nonzero control input, distributed
controllers are proposed in \cite{li2011trackingTAC,li2012adaptiveauto} to
solve the leader-follower consensus problem.
A common assumption in \cite{li2010consensus,ni2010leader,zhang2011optimal,you2011network,tuna2009conditions,seo2009consensus,ma2010necessary,li2011trackingTAC,li2012adaptiveauto} is that the dynamics of the agents
are identical and precisely known, which might be restrictive and not practical
in many circumstances.
In practical applications, the agents may be subject to certain parameter uncertainties
or unknown external disturbances.

This paper considers the distributed consensus problem of
multi-agent systems
with identical nominal linear dynamics
but subject to different matching uncertainties.
A typical example belonging to this scenario is a network of mass-spring
systems with different masses or unknown spring constants.
Due to the existence of the nonidentical uncertainties which may
be time-varying, nonlinear and unknown,
the multi-agent systems discussed in this paper are essentially
heterogeneous. The heterogeneous multi-agent systems in this paper
contain the homogeneous linear multi-agent systems studied
in \cite{li2010consensus,ni2010leader,zhang2011optimal,you2011network,tuna2009conditions,seo2009consensus,ma2010necessary}
as a special case where the uncertainties do not exist.
Note that because of the existence of the uncertainties,
the consensus problem in this case becomes quite challenging to solve
and the consensus algorithms given in \cite{li2010consensus,ni2010leader,zhang2011optimal,you2011network,tuna2009conditions,seo2009consensus,ma2010necessary}
are not applicable any more.

In this paper, we present a systematic procedure to
address the distributed robust consensus problem
of multi-agent systems with matching uncertainties
for both the cases
without and with a leader of possibly nonzero control input.
First, we consider the case where the communication graph is undirected and connected.
A distributed continuous
static consensus protocol based on the relative states
of neighboring agents is designed,
under which the consensus error is uniformly ultimately
bounded and exponentially converges to a
small residual set.
Note that
the design of this protocol relies on
the eigenvalues of the Laplacian matrix
and the upper bounds of the matching uncertainties.
In order to remove these requirements,
a fully distributed adaptive protocol is further designed,
under which the residual set of the consensus error is also given.
One desirable feature is that for both the static and adaptive protocols,
the residual sets of the consensus error
can be made to be reasonably small by properly selecting the design parameters
of the protocols and the convergence rates of the consensus
error are explicitly given.
Next, we extend to consider the case where there exists a leader with nonzero
control input. Here we study the general case
where the leader's control input
is not available to any follower,
which imposes additional difficulty.
Distributed static and adaptive consensus
protocols based on the relative state
information are proposed and designed to ensure
that the consensus error can converge
to residual sets which are explicitly given
and adjustable.
The case where the external disturbances
associated with the agent dynamics are bounded
and do not
satisfy the matching condition is also examined.
The proposed consensus protocols are redesigned
to guarantee the boundedness of the consensus
error.
The existence conditions of the
consensus protocols proposed in this paper are discussed.
It is pointed out that a sufficient condition
of the existence of the protocols
is that each agent is stabilizable.

It is worth mentioning that
in related works \cite{Das20102014,zhang2012adaptive},
the distributed tracking problem of multi-agent
systems with unknown nonlinear dynamics are discussed.
Compared to \cite{Das20102014,zhang2012adaptive},
the contribution of this paper is at least three-fold.
First, the agents in \cite{Das20102014,zhang2012adaptive}
are restricted to be first-order and special high-order systems.
It is far from trivial to extend \cite{Das20102014,zhang2012adaptive}
to solve the consensus problem
of the general high-order multi-agent systems with matching uncertainties
as in this paper. Second,
contrary to \cite{Das20102014,zhang2012adaptive} which consider only
the case with a leader, consensus for both the cases with and without a leader
is addressed in this paper. Third, the design of the protocols
in \cite{Das20102014,zhang2012adaptive} depends on global information
of the communication graph. In contrast, the adaptive consensus
protocols proposed in this paper are fully distributed,
which do not require any global information.


The rest of this paper is organized as follows. Some
useful results of graph theory
are reviewed in Section 2. The distributed robust leaderless consensus problem
is discussed in Section 3 for the case with an undirected graph.
The robust leader-follower consensus problem is addressed
in Section 4 for the case where there exists a leader with unknown control input.
The robustness of the proposed consensus protocols
with respect to external disturbances which do not satisfy the matching condition is
discussed in Section 5.
Simulation
examples are presented for illustration in Section 6.
Conclusions are drawn in
Section 7.

\section{Notation and Graph Theory}

$I_N$ represents the identity matrix of
dimension $N$.
Denote by $\mathbf{1}$ a column vector
with all entries equal to one. ${\rm{diag}}(A_1,\cdots,A_n)$
represents a block-diagonal matrix with matrices $A_i,i=1,\cdots,n,$
on its diagonal. 
$A\otimes B$ denotes the Kronecker product of matrices $A$ and $B$.
For a vector $x\in\mathbf{R}^n$, let
$\|x\|$ denote its 2-norm.
For a symmetric matrix $A$, $\lambda_{\min}(A)$ and $\lambda_{\max}(A)$
denote, respectively, the minimum and maximum eigenvalues
of $A$.

A directed graph $\mathcal {G}$ is a pair $(\mathcal {V}, \mathcal
{E})$, where $\mathcal {V}=\{v_1,\cdots,v_N\}$ is a nonempty finite
set of nodes and $\mathcal {E}\subseteq\mathcal {V}\times\mathcal
{V}$ is a set of edges, in which an edge is represented by an
ordered pair of distinct nodes. For an edge $(v_i,v_j)$, node $v_i$
is called the parent node, node $v_j$ the child node, and $v_i$ is a
neighbor of $v_j$. A graph with the property that
$(v_i,v_j)\in\mathcal {E}$ implies $(v_j, v_i)\in\mathcal {E}$ for
any $v_i,v_j\in\mathcal {V}$ is said to be undirected. A path from
node $v_{i_1}$ to node $v_{i_l}$ is a sequence of ordered edges of
the form $(v_{i_k}, v_{i_{k+1}})$, $k=1,\cdots,l-1$. A subgraph
$\mathcal {G}_s=(\mathcal {V}_s, \mathcal {E}_s)$ of $\mathcal {G}$
is a graph such that $\mathcal {V}_s\subseteq\mathcal {V}$ and
$\mathcal {E}_s\subseteq \mathcal {E}$.
A directed graph contains a directed spanning tree if
there exists a node called the root, which has no parent node, such
that the node has directed paths to all other nodes in the graph.

The adjacency matrix $\mathcal {A}=[a_{ij}]\in\mathbf{R}^{N\times
N}$ associated with the directed graph $\mathcal {G}$ is defined by
$a_{ii}=0$, $a_{ij}=1$ if $(v_j,v_i)\in\mathcal {E}$ and $a_{ij}=0$
otherwise. The Laplacian matrix $\mathcal {L}=[\mathcal
{L}_{ij}]\in\mathbf{R}^{N\times N}$ is defined as $\mathcal
{L}_{ii}=\sum_{j\neq i}a_{ij}$ and $\mathcal {L}_{ij}=-a_{ij}$,
$i\neq j$. For undirected graphs, both $\mathcal {A}$ and $\mathcal
{L}$ are symmetric.

{\bf Lemma 1} \cite{ren2005consensus}~
Zero is an eigenvalue of $\mathcal {L}$ with $\mathbf{1}$ as a
right eigenvector and all nonzero eigenvalues have positive real
parts. Furthermore, zero is a simple eigenvalue of $\mathcal {L}$ if
and only if $\mathcal {G}$ has a directed spanning tree.

\section{Distributed Robust Leaderless Consensus}

In this paper, we consider a network of $N$ autonomous agents
with identical nominal linear dynamics but
subject to heterogeneous uncertainties.
The dynamics of the $i$-th agent
are described by
\begin{equation}\label{1c1}
\dot{x}_i=Ax_i+Bu_i+H_i(x_i,t) + \nu_i(t), ~ i=1,\cdots,N,
\end{equation}
where $x_i\in\mathbf{R}^n$ is the state,
$u_i\in\mathbf{R}^{p}$ is the control input,
$A$ and $B$ are constant known matrices with
compatible dimensions, and
$H_i(x_i,t)\in\mathbf{R}^n$ and $\nu_i(t)\in\mathbf{R}^n$
denote, respectively, the parameter uncertainties
and external disturbances associated with the $i$-th agent,
which are assumed to
satisfy the following standard matching condition \cite{wheeler1998sliding,edwards1998sliding}.

{\bf Assumption 1}~There exist
functions $\hat{H}_i(x_i,t)$ and $\hat{\nu}_i(t)$ such
that $H_i(x_i,t)=B\hat{H}_i(x_i,t)$ and $\nu_i(t)=B\hat{\nu}_i(t)$, $i=1,\cdots,N$.

By letting $f_i(x_i,t)=\hat{H}_i(x_i,t)+\hat{\nu}_i$ represent the lumped uncertainty
of the $i$-th agent, \dref{1c1} can be rewritten into
\begin{equation}\label{1c}
\dot{x}_i=Ax_i+B[u_i+f_i(x_i,t)], ~ i=1,\cdots,N.
\end{equation}

In the previous related works
\cite{li2010consensus,ni2010leader,li2011dynamic,seo2009consensus,tuna2009conditions,zhang2011optimal,li2011adaptive,li2012adaptiveauto},
the agents are identical linear systems and free of uncertainties.
In contrast,
the agents \dref{1c} considered in this paper
are subject to nonidentical uncertainties,
which makes the resulting multi-agent systems
are essentially heterogeneous.
The agents \dref{1c} can recover the nominal linear agents
in \cite{li2010consensus,ni2010leader,li2011dynamic,seo2009consensus,tuna2009conditions,zhang2011optimal,li2011adaptive,li2012adaptiveauto}
when the uncertainties $f_i(x_i,t)$ do not exist.
Note that the existence of the uncertainties associated with the agents
makes the consensus problem quite challenging to solve,
as detailed in the sequel.

Regarding the bounds of the uncertainties
$f_i(x_i,t)$, we introduce the following assumption.

{\bf Assumption 2}~There exist continuous scalar valued
functions $\rho_i(x_i,t)$, $i=1,\cdots,N$, such
that $\|f_i(x_i,t)\|\leq\rho_i(x_i,t)$,
$i=1,\cdots,N$, for all $t\geq0$ and $x_i\in\mathbf{R}^n$.

The communication graph among the $N$ agents is represented by a
undirected graph $\mathcal {G}$, which is assumed
to be connected throughout this section.
The objective of this section is to solve the consensus problem
for the agents in \dref{1c}, i.e., to design distributed
consensus protocols such that
$\lim_{t\rightarrow \infty}\|x_i(t)- x_j(t)\|=0$, $
\forall\,i,j=1,\cdots,N$.

\subsection{Distributed Static Consensus Protocol}

Based on the relative states of neighboring agents,
the following distributed static consensus protocol is proposed:
\begin{equation}\label{ssd}
\begin{aligned}
u_i &=cK\sum_{j=0}^Na_{ij}(x_i-x_j)+\rho_i(x_i,t)\\
&\quad \times g (
K\sum_{j=0}^Na_{ij}(x_i-x_j) ),~i=1,\cdots,N,
\end{aligned}
\end{equation}
where
$c>0$ is the constant coupling gain,
$K\in\mathbf{R}^{p\times n}$ is the feedback gain matrix,
$a_{ij}$
is the $(i,j)$-th entry of the adjacency matrix $\mathcal {A}$
associated with $\mathcal {G}$, and the
nonlinear function $g(\cdot)$ is defined as follows:
for $w\in\mathbf{R}^n$,
\begin{equation}\label{satu}
g(w)=\begin{cases}\frac{w}{\|w\|} & \text{if}~ \rho_i(x_i,t)\|w\|>\kappa\\
\frac{w}{\kappa} & \text{if}~\rho_i(x_i,t)\|w\|\leq\kappa
\end{cases},
\end{equation}
where $\kappa$ is a small positive value.

Let $x=[x_1^T,\cdots,x_N^T]^T$ and
$\rho(x,t)={\rm{diag}}(\rho_1(x_1,t),\cdots,\\\rho_N(x_N,t))$. Using \dref{ssd} for
\dref{1c},
we can obtain the closed-loop network dynamics as
\begin{equation}\label{netss1}
\begin{aligned}
\dot{x}
&= (I_N\otimes A+c\mathcal {L}\otimes BK)x+(I_N\otimes B)F(x,t)\\
&\quad+[\rho(x,t)\otimes B]G(x),
\end{aligned}
\end{equation}
where $\mathcal {L}$ denotes the Laplacian matrix of $\mathcal {G}$, and
\begin{equation}\label{sa3}
F(x,t)\triangleq\begin{bmatrix}f_1(x_1,t)\\\vdots\\f_N(x_N,t)\end{bmatrix},
G(x)\triangleq\begin{bmatrix}g(K\sum_{j=1}^N \mathcal {L}_{1j}x_j)\\
\vdots\\g(K\sum_{j=1}^N \mathcal {L}_{Nj}x_j)\end{bmatrix}.
\end{equation}

Let $\xi=(M\otimes I_n)x$,
where $M=I_N-\frac{1}{N}\mathbf{1}\mathbf{1}^T$
and $\xi=[\xi_1^T,\cdots,\xi_N^T]^T$. It is easy
to see that $0$ is a simple eigenvalue of
$M$ with $\mathbf{1}$ as a
corresponding right eigenvector and 1 is the other eigenvalue with
multiplicity $N-1$. Then, it follows that $\xi=0$ if and only if
$x_1=\cdots=x_N$. Therefore, the consensus problem under the protocol
\dref{ssd} is solved if and only if
$\xi$ asymptotically converges to zero. Hereafter,
we refer to $\xi$ as the consensus error.
By noting that $\mathcal {L}M=\mathcal {L}$,
it is not difficult to obtain from \dref{netss1}
that the consensus error
$\xi$ satisfies
\begin{equation}\label{netss2}
\begin{aligned}
\dot{\xi}
&= (I_N\otimes A+c\mathcal {L}\otimes BK)\xi+(M\otimes B)F(x,t)\\
&\quad+[M\rho(x,t)\otimes B]G(\xi).
\end{aligned}
\end{equation}

The following result provides a sufficient condition
to design the consensus protocol \dref{ssd}.

{\bf Theorem 1}~
Suppose that the communication graph $\mathcal
{G}$ is undirected and connected
and Assumption 2 holds.
The parameters in the distributed protocol \dref{ssd}
are designed as $c\geq
\frac{1}{\lambda_2}$
and $K=-B^TP^{-1}$, where
$\lambda_2$ is the smallest nonzero eigenvalue of $\mathcal
{L}$ and $P>0$ is a solution to the following linear matrix
inequality (LMI):
\begin{equation}\label{alg1}
AP+PA^T-2BB^T<0,
\end{equation}
Then, the consensus error $\xi$ of \dref{netss2}
is uniformly ultimately bounded and exponentially
converges to the residual set
\begin{equation}\label{d}
\mathcal {D}_1\triangleq \{\xi : \|\xi\|^2\leq\frac{2\lambda_{\max}(P)N\kappa
}{\alpha\lambda_2}\},
\end{equation}
with a convergence rate faster than ${\rm{exp}(-\alpha t)}$, where
\begin{equation}\label{alf}
\alpha= \frac{-\lambda_{\max}(AP+PA^T-2 BB^T)}{\lambda_{\max}(P)}.
\end{equation}

{\bf Proof}~
Consider the following Lyapunov function candidate:
$$V_1=\frac{1}{2}\xi^T(\mathcal {L}\otimes P^{-1})\xi.$$
By the definition of $\xi$, it is easy to see that
$({\bf 1}^T\otimes I)\xi=0$.
For a connected graph $\mathcal {G}$,
it then follows from Lemma 1 that
\begin{equation}\label{lyas01}
V_1(\xi)\geq
\frac{1}{2}\lambda_2\xi^T(I_N\otimes P^{-1})\xi
\geq \frac{\lambda_2}{2\lambda_{\max}(P)}\|\xi\|^2.
\end{equation}
The time derivative of $V_1$
along the trajectory of \dref{netss1} is given by
\begin{equation}\label{lyas2}
\begin{aligned}
\dot{V}_1 
&=\xi^T(\mathcal {L}\otimes P^{-1}A+c\mathcal {L}^2\otimes P^{-1}BK)\xi\\
&\quad+\xi^T(\mathcal {L}\otimes P^{-1}B)F(x,t)\\
&\quad+\xi^T[\mathcal {L}\rho(x,t)\otimes P^{-1}B]G(\xi).
\end{aligned}
\end{equation}

By using Assumption 2, we can obtain that
\begin{equation}\label{lyas5}
\begin{aligned}
&\xi^T(\mathcal {L}\otimes P^{-1}B) F(x,t)\\
&\qquad\leq
\sum_{i=1}^N \|B^TP^{-1}\sum_{j=1}^N\mathcal {L}_{ij}\xi_j\| \| f_i(x_i,t)\|\\
&\qquad\leq
 \sum_{i=1}^N\rho_i(x_i,t) \|B^TP^{-1}\sum_{j=1}^N\mathcal {L}_{ij}\xi_j\|.
\end{aligned}
\end{equation}

Next, consider the following three cases.

i) $\rho_i(x_i,t)\|K\sum_{j=1}^N \mathcal {L}_{ij}\xi_j\|>\kappa$,
$i=1,\cdots,N$.

In this case, it follows from \dref{satu} and \dref{sa3} that
\begin{equation}\label{lyas4}
\begin{aligned}
&\xi^T[\mathcal {L}\rho(x,t)\otimes P^{-1}B]G(\xi)\\
&\qquad=-\sum_{i=1}^N \rho_i(x_i,t)\|B^TP^{-1}\sum_{j=1}^N\mathcal {L}_{ij}\xi_j\|.
\end{aligned}
\end{equation}
Substituting \dref{lyas4} and and \dref{lyas5} into \dref{lyas2} yields
$ \dot{V}_1 \leq
\frac{1}{2}\xi^T\mathcal {X}\xi$,
where
$\mathcal {X}=\mathcal {L}\otimes (P^{-1}A+A^TP^{-1})-2c\mathcal {L}^2\otimes P^{-1}BB^TP^{-1}$.

ii) $\rho_i(x_i,t)\|K\sum_{j=1}^N \mathcal {L}_{ij}\xi_j\|\leq\kappa$, $i=1,\cdots,N$.

In this case, we can get from \dref{satu} and \dref{sa3} that
\begin{equation}\label{lyasc8}
\begin{aligned}
&\xi^T[\mathcal {L}\rho(x,t)\otimes P^{-1}B]G(\xi)\\
&\qquad=-\sum_{i=1}^N \frac{\rho_i(x_i,t)}{\kappa}\|B^TP^{-1}\sum_{j=1}^N \mathcal {L}_{ij}\xi_j\|^2
\leq 0.
\end{aligned}
\end{equation}
Substituting \dref{lyas4}, \dref{lyas5}, and \dref{lyasc8} into \dref{lyas2} gives
\begin{equation}\label{lyasc9}
\begin{aligned}
\dot{V}_1 \leq \frac{1}{2}\xi^T\mathcal {X}\xi+N\kappa.
\end{aligned}
\end{equation}

iii) $\xi$ satisfies neither case i) nor case ii).

Without loss of generality, assume that
$\rho_i(x_i,t)\|K\sum_{j=1}^N \mathcal {L}_{ij}\xi_j\|>\kappa$, $i=1,\cdots,l$, and
$\rho_i(x_i,t)\|K\sum_{j=1}^N \mathcal {L}_{ij}\xi_j\|\leq\kappa$, $i=l+1,\cdots,N$,
where $2\leq l\leq N-1$.
By combing \dref{lyas4} and \dref{lyasc8}, in this case we can get that
\begin{equation}\label{lyasc10}
\begin{aligned}
&\xi^T[\mathcal {L}\rho(x,t)\otimes P^{-1}B]G(\xi)\\
&\qquad\leq -\sum_{i=1}^l \rho_i(x_i,t)\|B^TP^{-1}\sum_{j=1}^N\mathcal {L}_{ij}\xi_j\|.
\end{aligned}
\end{equation}
Then, it follows from \dref{lyas2}, \dref{lyas4}, \dref{lyasc10}, and \dref{lyas5} that
$
\dot{V}_1 \leq \frac{1}{2}\xi^T\mathcal {X}\xi+(N-l)\kappa.
$

Therefore, by analyzing the above three cases, we get that
$\dot{V}_1$ satisfies \dref{lyasc9} for all $\xi\in\mathbf{R}^{Nn}$.
Note that \dref{lyasc9} can be rewritten as
\begin{equation}\label{lyasc11}
\begin{aligned}
\dot{V}_1 &\leq -\alpha V_1+\alpha V_1+\frac{1}{2}\xi^T\mathcal {X}\xi+N\kappa\\
&=-\alpha V_1+\frac{1}{2}\xi^T(\mathcal {X}+\alpha\mathcal {L}\otimes P^{-1})\xi+N\kappa,
\end{aligned}
\end{equation}
where $\alpha>0$.

Because $\mathcal {G}$ is connected, it follows from Lemma 1 that
zero is a simple eigenvalue of
$\mathcal {L}$ and all the other eigenvalues are positive.
Let $U=\left[\begin{smallmatrix}
\frac{\mathbf{1}}{\sqrt{N}} & Y_1
\end{smallmatrix}\right]$ and $U^T=\left[\begin{smallmatrix}
\frac{\mathbf{1}^T}{\sqrt{N}} \\ Y_2
\end{smallmatrix}\right]$, with $Y_1\in\mathbf{R}^{N\times(N-1)}$, $Y_2\in\mathbf{R}^{(N-1)\times N}$,
be such unitary matrices that
$U^{T}\mathcal {L}U=\Lambda\triangleq{\rm{diag}}(0,\lambda_2,\cdots,\lambda_N)$,
where $\lambda_2\leq\cdots\leq\lambda_N$ are the nonzero eigenvalues of $\mathcal {L}$.
Let $\bar{\xi}\triangleq[\bar{\xi}_1^T,\cdots,\bar{\xi}_N^T]^T=(U^T\otimes P^{-1})\xi$.
By the definitions of $\xi$ and $\bar{\xi}$, it is easy to see that
$
\bar{\xi}_1=(\frac{\mathbf{1}^T}{\sqrt{N}}\otimes P^{-1})\xi
=(\frac{\mathbf{1}^T}{\sqrt{N}}M\otimes P^{-1})x=0.
$
Then, it follows that
\begin{equation}\label{lyasc12}
\begin{aligned}
&\xi^T(\mathcal {X}+\alpha\mathcal {L}\otimes P^{-1})\xi\\
&\qquad=\sum_{i=2}^{N}\lambda_i\bar{\xi}_i^T(AP+PA^T+\alpha P-2c\lambda_i BB^T)\bar{\xi}_i\\
&\qquad\leq \sum_{i=2}^{N}\lambda_i\bar{\xi}_i^T(AP+PA^T+\alpha P-2BB^T)\bar{\xi}_i.
\end{aligned}
\end{equation}
Because $\alpha= \frac{-\lambda_{\max}(AP+PA^T-2 BB^T)}{\lambda_{\max}(P)}$,
we can see from \dref{lyasc12} that $\xi^T(\mathcal {X}+\alpha\mathcal {L}\otimes P^{-1})\xi\leq0$.
Then, we can get from \dref{lyasc11}
that
\begin{equation}\label{lyasc13}
\begin{aligned}
\dot{V}_1
\leq-\alpha V_1+N\kappa.
\end{aligned}
\end{equation}
By using the well-known Comparison lemma (Lemma 3.4 in \cite{khalil2002nonlinear}),
we can obtain from \dref{lyasc13} that
\begin{equation}\label{lyasc14}
\begin{aligned}
V_1(\xi)
\leq [V_1(\xi(0))-\frac{N\kappa}{\alpha}]{\rm{exp}(-\alpha t)}+\frac{N\kappa}{\alpha},
\end{aligned}
\end{equation}
which, by \dref{lyas01}, implies
that $\xi$ exponentially converges to the residual set
$\mathcal {D}_1$ in \dref{d} with a convergence rate not less than
${\rm{exp}(-\alpha t)}$.
\hfill $\blacksquare$

{\bf Remark 1}~
The distributed consensus protocol \dref{ssd} consists of
a linear part and a nonlinear part,
where the term $\rho_i(x_i,t)g(
K\sum_{j=1}^Na_{ij}(x_i-x_j))$
is used to suppress the effect of the uncertainties $f_i(x_i,t)$.
For the case where $f_i(x_i,t)=0$, we can accordingly remove
$\rho_i(x_i,t)g(
K\sum_{j=1}^Na_{ij}(x_i-x_j))$ from \dref{ssd},
which can recover the static consensus protocols
as in \cite{li2010consensus,li2011dynamic,zhang2011optimal}.
As shown in Proposition 2 of \cite{li2010consensus}, a necessary and
sufficient condition for the existence of a $P>0$ to the LMI
\dref{alg1} is that $(A,B)$ is stabilizable. Therefore, a sufficient
condition for the existence of \dref{ssd} satisfying Theorem 1 is
that $(A,B)$ is stabilizable.
Note that in Theorem 1 the parameters $c$ and $K$ of
\dref{ssd} are independently designed.

Note that the nonlinear component $g(\cdot)$
in \dref{satu} is continuous,
which is actually a continuous approximation,
via the boundary layer concept \cite{edwards1998sliding,khalil2002nonlinear},
of the discontinuous
function
$\hat{g}(w)=\begin{cases}\frac{w}{\|w\|} & \text{if}~\|w\|\neq 0\\
0 & \text{if}~\|w\|=0
\end{cases}.$
The value of $\kappa$ in \dref{satu} defines the size of
the boundary layer. As $\kappa\rightarrow 0$,
the continuous function $g(\cdot)$
approaches the discontinuous function $\hat{g}(\cdot)$.


{\bf Corollary 1}~
Assume that $\mathcal {G}$ is connected and Assumption 2 holds. The consensus
error $\xi$ converges to zero under the discontinuous
consensus protocol:
\begin{equation}\label{sscdis}
\begin{aligned}
u_i &=cK\sum_{j=1}^Na_{ij}(x_i-x_j)+
\rho_i(x_i,t)\\
&\quad\times\hat{g}(K\sum_{j=1}^Na_{ij}(x_i-x_j) ),~i=1,\cdots,N,
\end{aligned}
\end{equation}
where $c$ and $K$ are chosen as in Theorem 1.

{\bf Remark 2}~
An inherent drawback of the discontinuous
protocol \dref{sscdis}
is that it will result in
the undesirable chattering effect
in real implementation,
due to imperfections in switching devices
\cite{young1999control,edwards1998sliding}.
The effect of chattering is avoided
by using the continuous protocol \dref{ssd}.
The cast is that the protocol \dref{ssd} does no
guarantee asymptotic stability but rather uniform
ultimate boundedness of the consensus error $\xi$.
Note that the residual set $\mathcal {D}_1$
of $\xi$
depends on the smallest nonzero eigenvalue of
$\mathcal {L}$, the number of agents,
the largest eigenvalue of $P$,
and the size $\kappa$ of the boundary layer.
By choosing a sufficiently small $\kappa$,
the consensus error $\xi$ under the
protocol \dref{ssd} can
converge to an arbitrarily small neighborhood of zero,
which is acceptable in most applications.

\subsection{Distributed Adaptive Consensus Protocol}

In the last subsection, the design of the distributed protocol
\dref{ssd} relies on the minimal nonzero eigenvalue $\lambda_2$ of
$\mathcal {L}$ and the upper bounds $\rho_i(x_i,t)$ of
the matching uncertainties $f_i(x_i,t)$.
However,
$\lambda_2$ is global information in the sense that
each agent has to know the entire
communication graph to compute it.
Besides, the bounds $\rho_i(x_i,t)$
of the uncertainties $f_i(x_i,t)$
might not be easily obtained in some cases,
e.g., $f_i(x_i,t)$
contains certain unknown
external disturbances.
In this subsection, we will implement some adaptive
control ideas to compensate the lack of
$\lambda_2$ and $\rho_i(x_i,t)$
and thereby to solve the consensus problem using only
the local information available to each agent.

Before moving forward, we introduce a modified
assumption regarding the bounds of the lumped uncertainties
$f_i(x_i,t)$, $i=1,\cdots,N$.

{\bf Assumption 3}~ There are positive constants $d_i$ and
$e_i$ such that $\|f_i(x_i,t)\|\leq d_i+e_i\|x_i\|$, $i=1,\cdots,N$.

Based on the local state information of neighboring agents,
we propose the following distributed adaptive
protocol to each agent:
\begin{equation}\label{ssca}
\begin{aligned}
u_i &=\bar{d}_iK\sum_{j=1}^Na_{ij}(x_i-x_j)+r (
K\sum_{j=1}^Na_{ij}(x_i-x_j)),\\
\dot{\bar{d}}_i &= \tau_i [-\varphi_i\bar{d}_i+(\sum_{j=1}^Na_{ij}(x_i-x_j)^T)\Gamma
(\sum_{j=1}^Na_{ij}(x_i-x_j))\\
&\quad+\|K\sum_{j=1}^Na_{ij}(x_i-x_j)\|],\\
\dot{\bar{e}}_i &=\epsilon_i [-\psi_i\bar{e}_i+\|K\sum_{j=1}^Na_{ij}(x_i-x_j)\|\|x_i\|],
i=1,\cdots,N,
\end{aligned}
\end{equation}
where
$\bar{d}_i(t)$ and $\bar{e}_i(t)$
are the adaptive gains associated with the $i$-th agent,
$\Gamma\in\mathbf{R}^{n\times n}$ is the feedback gain
matrix, $\tau_i$ and $\epsilon_i$ are positive scalars,
$\varphi_i$ and $\psi_i$ are small positive constants chosen
by the designer, the nonlinear function $r(\cdot)$
is defined as follows: for $w\in\mathbf{R}^n$,
\begin{equation}\label{satua}
r(w)=\begin{cases} \frac{w(\bar{d}_i+\bar{e}_i\|x_i\|)}{\|w\|} & \text{if}~(\bar{d}_i+\bar{e}_i\|x_i\|)\|w\|>\kappa \\
\frac{w(\bar{d}_i+\bar{e}_i\|x_i\|)^2}{\kappa} & \text{if}~(\bar{d}_i+\bar{e}_i\|x_i\|)\|w\|\leq\kappa
\end{cases},
\end{equation}
and the rest of the variables
are defined as in \dref{ssd}.

Let the consensus error $\xi$ be defined as in \dref{netss2} and $\overline{D}={\rm{diag}}(\bar{d}_1,\cdots,\bar{d}_N)$.
Then, it is not difficult to get from \dref{1c} and \dref{ssca}
that the closed-loop network dynamics can be written as
\begin{equation}\label{netsca1}
\begin{aligned}
\dot{\xi}
&= (I_N\otimes A+M\overline{D}\mathcal {L}\otimes BK)\xi+(M\otimes B)F(x,t)\\
&\quad
+(M\otimes B)R(\xi),\\
\dot{\bar{d}}_i &= \tau_i [-\varphi_i\bar{d}_i+(\sum_{j=1}^N\mathcal {L}_{ij}\xi_j^T)\Gamma
(\sum_{j=1}^N\mathcal {L}_{ij}\xi_j)+\|K\sum_{j=1}^N\mathcal {L}_{ij}\xi_j\| ],\\
\dot{\bar{e}}_i &=\epsilon_i [-\psi_i\bar{e}_i+\|K\sum_{j=1}^N\mathcal {L}_{ij}\xi_j\|\|x_i\| ],
~i=1,\cdots,N,
\end{aligned}
\end{equation}
where
\begin{equation}\label{sa4}
R(\xi)\triangleq\begin{bmatrix}r(K\sum_{j=1}^N \mathcal {L}_{1j}\xi_j)\\
\vdots\\r(K\sum_{j=1}^N \mathcal {L}_{Nj}\xi_j)\end{bmatrix},
\end{equation}
and the rest of the variables are defined as in \dref{netss1}.

To establish the ultimate boundedness of
the states $\xi$, $\bar{d}_i$, and $\bar{e}_i$
of \dref{netsca1}, we use the following Lyapunov function candidate
\begin{equation}\label{lyaas1}
V_2=\frac{1}{2}\xi^T(\mathcal {L}\otimes P^{-1})\xi+\sum_{i=1}^N\frac{\tilde{d}_i^2}{2\tau_i}
+\sum_{i=1}^N\frac{\tilde{e}_i^2}{2\epsilon_i},
\end{equation}
where $\tilde{e}_i=\bar{e}_i-e_i$,
$\tilde{d}_i=\bar{d}_i-\beta$, $i=1,\cdots,N$, and
$\beta\geq \max_{i=1,\cdots,N} \{d_i,\frac{1}{\lambda_2}\}$.


{\bf Theorem 3}~Suppose that $\mathcal{G}$ is connected and
Assumption 3 holds.
The feedback gain matrices of the distributed
adaptive protocol \dref{ssca}
are designed as
$K=-B^TP^{-1}$ and $\Gamma=P^{-1}BB^TP^{-1}$, where
$P>0$ is a solution to the LMI \dref{alg1}.
Then, both the consensus error $\xi$
and the adaptive gains $\bar{d}_i$
and $\bar{e}_i$, $i=1,\cdots,N$, in \dref{netsca1}
are uniformly ultimately bounded and the following statements hold.

\begin{itemize}
\item[i)]
For any $\varphi_i$ and $\psi_i$,
$\xi$, $\tilde{d}_i$,
and $\tilde{e}_i$ exponentially converge to the residual set
\begin{equation}\label{d2}
\mathcal {D}_2\triangleq \{\xi,\tilde{d}_i,\tilde{e}_i:V_2<
\frac{1}{2\delta}\sum_{i=1}^N(\beta^2\varphi_i+e_i^2\psi_i)+\frac{N\kappa}{4\delta}\},
\end{equation}
with a convergence rate faster than ${\rm{exp}(-\delta t)}$,
where $ \delta\triangleq \min_{i=1,\cdots,N}\{\alpha,\varphi_i \tau_i,\psi_i\epsilon_i\}$
and $\alpha$ is defined as in \dref{alf}.
\item[ii)]
If small $\varphi_i$ and $\psi_i$ satisfy
$\varrho\triangleq\max_{i=1,\cdots,N}\{\varphi_i \tau_i,\psi_i\epsilon_i\}<\alpha$,
then
in addition to i), $\xi$ exponentially converges to the residual set
\begin{equation}\label{d2k}
\mathcal {D}_3\triangleq \{\xi:\|\xi\|^2\leq
\frac{\lambda_{\max}(P)}{\lambda_2(\alpha-\varrho)}[\sum_{i=1}^N(\beta^2\varphi_i+e_i^2\psi_i)+\frac{1}{2}N\kappa]\}.
\end{equation}
\end{itemize}
with a convergence rate faster than ${\rm{exp}(-\varrho t)}$.

{\bf Proof}~
The time derivative of $V_2$ along \dref{netsca1}
can be obtained as
\begin{equation}\label{lyasac1}
\begin{aligned}
\dot{V}_2
&=\xi^T[(\mathcal {L}\otimes P^{-1}A+\mathcal {L}\widetilde{D}\mathcal {L}\otimes P^{-1}BK)\xi
\\&\quad+(\mathcal {L}\otimes P^{-1}B)F(x,t)+(\mathcal {L}\otimes P^{-1}B)R(\xi)]\\
&\quad +\sum_{i=1}^N\tilde{d}_i [-\varphi_i(\tilde{d}_i+\beta)+(\sum_{j=1}^N\mathcal {L}_{ij}\xi_j^T)\Gamma(\sum_{j=1}^N\mathcal {L}_{ij}\xi_j)\\&\quad+\|K\sum_{j=1}^N\mathcal {L}_{ij}\xi_j\|]+\sum_{i=1}^N\tilde{e}_i [-\psi_i(\tilde{e}_i+e_i)
\\
&\quad +\|K\sum_{j=1}^N\mathcal {L}_{ij}\xi_j\|\|x_i\|],
\end{aligned}
\end{equation}
where $\widetilde{D}(t)={\rm{diag}}(\tilde{d}_1+\beta,\cdots,\tilde{d}_N+\beta)$.

By noting that $K=-BP^{-1}$, it is easy to get that
\begin{equation}\label{lyasac2}
\begin{aligned}
&\xi^T(\mathcal {L}\widetilde{D}\mathcal {L}\otimes P^{-1}BK)\xi
\\&=
-\sum_{i=1}^N(\tilde{d}_i+\beta)(\sum_{j=1}^N\mathcal {L}_{ij}\xi_j)^TP^{-1}BB^TP^{-1}(\sum_{j=1}^N\mathcal {L}_{ij}\xi_j).
\end{aligned}
\end{equation}
In light of Assumption 3, we can obtain that
\begin{equation}\label{lyasac22}
\begin{aligned}
&\xi^T(\mathcal {L}\otimes P^{-1}B)F(x,t)
\\&\qquad\leq\sum_{j=1}^N(d_i+e_i\|x_i\|)\|B^TP^{-1}\sum_{j=1}^N \mathcal {L}_{ij}\xi_j\|.
\end{aligned}
\end{equation}

In what follows, we consider three cases.

i) $(\bar{d}_i+\bar{e}_i\|x_i\|)\|K\sum_{j=1}^N\mathcal {L}_{ij}\xi_j\|>\kappa$, $i=1,\cdots,N$.

In this case, we can get from \dref{satua} and \dref{sa4} that
\begin{equation}\label{lyasac3}
\begin{aligned}
&\xi^T(\mathcal {L}\otimes P^{-1}B)R(\xi)
\\&=-\sum_{i=1}^N[\tilde{d}_i+\beta+(\tilde{e}_i+e_i)\|x_i\|]\|B^TP^{-1}\sum_{j=1}^N \mathcal {L}_{ij}\xi_j\|.
\end{aligned}
\end{equation}
Substituting \dref{lyasac2}, \dref{lyasac22},
and \dref{lyasac3} into \dref{lyasac1} yields
$$\begin{aligned}
\dot{V}_2
&\leq \frac{1}{2}\xi^T\mathcal {Y}\xi-\sum_{i=1}^N(\beta-d_i)\|B^TP^{-1}\sum_{j=1}^N \mathcal {L}_{ij}\xi_j\|\\
&\quad-\frac{1}{2}\sum_{i=1}^N(\varphi_i\tilde{d}_i^2+\psi_i\tilde{e}_i^2)
+\frac{1}{2}\sum_{i=1}^N(\beta^2\varphi_i+e_i^2\psi_i)\\
&\leq \frac{1}{2}\xi^T\mathcal {Y}\xi-\frac{1}{2}\sum_{i=1}^N(\varphi_i\tilde{d}_i^2+\psi_i\tilde{e}_i^2)
+\frac{1}{2}\sum_{i=1}^N(\beta^2\varphi_i+e_i^2\psi_i),
\end{aligned}$$
where
$\mathcal {Y}\triangleq\mathcal {L}\otimes (P^{-1}A+A^TP^{-1})-2\beta \mathcal {L}^2\otimes P^{-1}BB^TP^{-1}$
and
we have used the facts that $\beta\geq \max_{i=1,\cdots,N} d_i $
and $-\tilde{d}_i^2-\tilde{d}_i\beta\leq-\frac{1}{2}\tilde{d}_i^2+\frac{1}{2}\beta^2$.

ii) $(\bar{d}_i+\bar{e}_i\|x_i\|)\|K\sum_{j=1}^N\mathcal {L}_{ij}\xi_j\|\leq\kappa$, $i=1,\cdots,N$.

In this case, we can get from \dref{satua} and \dref{sa4} that
\begin{equation}\label{lyasac4}
\begin{aligned}
&\xi^T(\mathcal {L}\otimes P^{-1}B)R(\xi)
\\&\quad=-\sum_{i=1}^N\frac{(\bar{d}_i+\bar{e}_i\|x_i\|)^2}{\kappa}\|B^TP^{-1}\sum_{j=1}^N \mathcal {L}_{ij}\xi_j\|^2.
\end{aligned}
\end{equation}
Then, it follows from \dref{lyasac2}, \dref{lyasac22},
\dref{lyasac4}, and \dref{lyasac1} that
\begin{equation}\label{lyasac5}
\begin{aligned}
\dot{V}_2
&\leq \frac{1}{2}\xi^T\mathcal {Y}\xi
-\frac{1}{2}\sum_{i=1}^N(\varphi_i\tilde{d}_i^2+\psi_i\tilde{e}_i^2)
\\&\quad+\frac{1}{2}\sum_{i=1}^N(\beta^2\varphi_i+e_i^2\psi_i)+\frac{1}{4}N\kappa,
\end{aligned}
\end{equation}
where we have used the fact that
$
-\frac{(\bar{d}_i+\bar{e}_i\|x_i\|)^2}{\kappa}\|B^TP^{-1}\\\sum_{j=1}^N \mathcal {L}_{ij}\xi_j\|^2
+(\bar{d}_i+\bar{e}_i\|x_i\|)\|B^TP^{-1}\sum_{j=1}^N \mathcal {L}_{ij}\xi_j\|\leq \frac{1}{4}\kappa,
$
for $(\bar{d}_i+\bar{e}_i\|x_i\|)\|K\sum_{j=1}^N\mathcal {L}_{ij}\xi_j\|\leq\kappa$, $i=1,\cdots,N$.

iii) $(\bar{d}_i+\bar{e}_i\|x_i\|)\|K\sum_{j=1}^N\mathcal {L}_{ij}\xi_j\|>\kappa$, $i=1,\cdots,l$, and
$(\bar{d}_i+\bar{e}_i\|x_i\|)\|K\sum_{j=1}^N\mathcal {L}_{ij}\xi_j\|\leq\kappa$, $i=l+1,\cdots,N$,
where $2\leq l\leq N-1$.

By following similar
steps in the two cases above, it is not difficult to get that
$$\begin{aligned}
\dot{V}_2
&\leq\frac{1}{2}\xi^T\mathcal {Y}\xi
-\frac{1}{2}\sum_{i=1}^N(\varphi_i\tilde{d}_i^2+\psi_i\tilde{e}_i^2)
\\&\qquad+\frac{1}{2}\sum_{i=1}^N(\beta^2\varphi_i+e_i^2\psi_i)
+\frac{1}{4}(N-l)\kappa.
\end{aligned}$$

Therefore, based on the above three cases, we can get that
$\dot{V}_2$ satisfies \dref{lyasac5} for all $\xi\in\mathbf{R}^{Nn}$.
Note that \dref{lyasac5} can be rewritten into
\begin{equation}\label{lyasac61}
\begin{aligned}
\dot{V}_2
&\leq -\delta V_2+\delta V_2+\frac{1}{2}\xi^T\mathcal {Y}\xi-\frac{1}{2}\sum_{i=1}^N(\varphi_i\tilde{d}_i^2+\psi_i\tilde{e}_i^2)
\\&\quad+\frac{1}{2}\sum_{i=1}^N(\beta^2\varphi_i+e_i^2\psi_i)+\frac{1}{4}N\kappa\\
&=-\delta V_2+\frac{1}{2}\xi^T(\mathcal {Y}+\delta \mathcal {L}\otimes P^{-1})\xi
-\frac{1}{2}\sum_{i=1}^N[(\varphi_i-\frac{\delta}{\tau_i})\tilde{d}_i^2\\
&\quad+(\psi_i-\frac{\delta}{\epsilon_i})\tilde{e}_i^2)]+\frac{1}{2}\sum_{i=1}^N(\beta^2\varphi_i+e_i^2\psi_i)+\frac{1}{4}N\kappa.
\end{aligned}
\end{equation}
Because $\beta\lambda_2\geq1$ and $0<\delta\leq\alpha$,
by following similar steps in the proof of Theorem 1,
we can show that $\xi^T(\mathcal {Y}+\delta \mathcal {L}\otimes P^{-1})\xi\leq0$. Further, by noting that
$\delta\leq\min_{i=1,\cdots,N}\{\varphi_i \tau_i,\psi_i\epsilon_i\}$,
it follows from \dref{lyasac61} that
\begin{equation}\label{lyasac62}
\begin{aligned}
\dot{V}_2
&\leq
-\delta V_2+\frac{1}{2}\sum_{i=1}^N(\beta^2\varphi_i+e_i^2\psi_i)+\frac{1}{4}N\kappa,
\end{aligned}
\end{equation}
which implies that
\begin{equation}\label{lyasc142}
\begin{aligned}
V_2
&\leq [V_2(0)-\frac{N\kappa}{4\delta}-
\frac{1}{2\delta}\sum_{i=1}^N(\beta^2\varphi_i+e_i^2\psi_i)]{\rm{exp}(-\delta t)}\\
&\quad+
\frac{1}{2\delta}\sum_{i=1}^N(\beta^2\varphi_i+e_i^2\psi_i)+\frac{N\kappa}{4\delta}.
\end{aligned}
\end{equation}
Therefore, $V_2$ exponentially converges to the residual set
$\mathcal {D}_2$ in \dref{d2} with a convergence rate faster than ${\rm{exp}(-\delta t)}$,
which, in light of $V_2\geq \frac{\lambda_2}{2\lambda_{\max}(P)}\|\xi\|^2$,
implies that $\xi$,
$\bar{d}_i$, and $\bar{e}_i$ are uniformly ultimately bounded.

Next, if $\varrho\triangleq\max_{i=1,\cdots,N}\{\varphi_i \tau_i,\psi_i\epsilon_i\}\leq\alpha$,
then we can choose $\delta=\varrho$
and it is easy to see that $\mathcal {D}_3$ increase as $\varphi_i$ and $\psi_i$ decrease.
However, for the case where
$\varrho<\alpha$, we can obtain a
smaller residual set for $\xi$ by rewriting \dref{lyasac61} into
\begin{equation}\label{lyasac71}
\begin{aligned}
\dot{V}_2
&\leq
-\varrho V_2+\frac{1}{2}\xi^T(\mathcal {Y}+\alpha\mathcal {L}\otimes P^{-1})\xi+\frac{1}{4}N\kappa\\
&\quad-\frac{\alpha-\varrho}{2}\xi^T(\mathcal {L}\otimes P^{-1})\xi
+\frac{1}{2}\sum_{i=1}^N(\beta^2\varphi_i+e_i^2\psi_i)\\
&\leq -\varrho V_2-\frac{\lambda_2(\alpha-\varrho)}{2\lambda_{\max}(P)}\|\xi\|^2
+\frac{1}{2}\sum_{i=1}^N(\beta^2\varphi_i+e_i^2\psi_i)\\&\quad+\frac{1}{4}N\kappa.
\end{aligned}
\end{equation}
Obviously, it follows from \dref{lyasac71} that $\dot{V}_2\leq -\varrho V_2$ if
$\|\xi\|^2>\frac{\lambda_{\max}(P)}{\lambda_2(\alpha-\varrho)}[\sum_{i=1}^N(\beta^2\varphi_i+e_i^2\psi_i)+\frac{1}{2}N\kappa].$
Then, by noting $V_2\geq \frac{\lambda_2}{2\lambda_{\max}(P)}\|\xi\|^2$,
we can get that if $\varrho\leq\alpha$ then
$\xi$ exponentially converges to the residual set $\mathcal {D}_3$ in \dref{d2k}
with a convergence rate faster than ${\rm{exp}(-\varrho t)}$.
\hfill $\blacksquare$

{\bf Remark 3}~
It is worth mentioning that
adding $-\varphi_i\bar{d}_i$ and
$-\psi_i\bar{e}_i$ into \dref{ssca}
is essentially motivated by the so-called $\sigma$-modification
technique in \cite{ioannou1984instability,wheeler1998sliding},
which plays a vital role to guarantee the ultimate
boundedness of the consensus error $\xi$
and the adaptive gains $\bar{d}_i$ and $\bar{e}_i$.
From \dref{d2} and \dref{d2k}, we can observe that
the residual sets $\mathcal {D}_2$ and $\mathcal {D}_3$
decrease as $\kappa$ decreases.
Given $\kappa$, smaller $\varphi_i$ and $\psi_i$ give a smaller bound for
$\xi$ and at the same time yield a larger bound for
$\bar{d}_i$ and $\bar{e}_i$. For the case
where $\varphi_i=0$ and $\psi_i=0$, $\bar{d}_i$ and $\bar{e}_i$
will tend to infinity. In real implementations, if large
$\bar{d}_i$ and $\bar{e}_i$ are acceptable, we can choose $\varphi_i$, $\psi_i$,
and $\kappa$ to be relatively small in order to guarantee
a small consensus error $\xi$.

{\bf Remark 4}~
Contrary to the static
protocol \dref{ssd}, the design of the
adaptive protocol \dref{ssca}
relies on only the agent dynamics,
requiring neither the minimal nonzero
eigenvalue of $\mathcal {L}$ nor the upper bounds of
of the uncertainties $f_i(x_i,t)$.
Thus, the adaptive controller \dref{ssca}
can be implemented by each agent in a fully distributed
fashion without requiring
any global information.

{\bf Remark 5}~
A special case of the uncertainties $f_i(x_i,t)$ satisfying Assumptions 1 and 2
is that there exist positive constants $d_i$ such
that $\|f_i(x_i,t)\|\leq d_i$, $i=1,\cdots,N$. For this case, the proposed protocols
\dref{ssd} and \dref{ssca} can be accordingly simplified.
For \dref{ssd}, we can simply replace $\rho_i(x_i,t)$ by $d_i$.
The adaptive protocol \dref{ssca} in this case can be modified into
\begin{equation}\label{sscas}
\begin{aligned}
u_i &=\bar{d}_iK\sum_{j=1}^Na_{ij}(x_i-x_j)+\bar{r} (
K\sum_{j=1}^Na_{ij}(x_i-x_j)),\\
\dot{\bar{d}}_i &= \tau_i [-\varphi_i\bar{d}_i+(\sum_{j=1}^Na_{ij}(x_i-x_j)^T)\Gamma
(\sum_{j=1}^Na_{ij}(x_i-x_j))\\&\quad+\|K\sum_{j=1}^Na_{ij}(x_i-x_j)\|],
~i=1,\cdots,N,
\end{aligned}
\end{equation}
where
the nonlinear function $\bar{r}(\cdot)$
is defined such that
$\bar{r}(w)=\begin{cases} \frac{w\bar{d}_i}{\|w\|}, & \text{if}~ \bar{d}_i\|w\|>\kappa \\
\frac{w \bar{d}_i}{\kappa}, & \text{if}~ \bar{d}_i\|w\|\leq\kappa
\end{cases}$
and the rest of the variables
are defined as in \dref{ssca}. 

\section{Distributed Robust Leader-Follower Consensus with a Leader of Nonzero
Control Input}

For the leaderless consensus problem in the previous section
where the communication graph is undirected,
the final consensus values reached by the agents under the protocols
\dref{ssd} and \dref{ssca} are generally difficult to be explicitly obtained.
The main difficulty lies in that the agents are subject
to uncertainties and the protocols \dref{ssd} and \dref{ssca}
are essentially nonlinear.
In this section, we consider the leader-follower consensus problem,
for which case the agents' states
are required to converge onto a reference trajectory.

Consider a network of $N+1$ agents consisting of
$N$ followers and one leader.
Without loss of generality, let the agent indexed by $0$
be the leader and the agents indexed by $1,\cdots,N$, be the
followers. The dynamics of the followers are described by
\dref{1c}. For simplicity, we assume that the leader has the nominal
linear dynamics, given by
\begin{equation}\label{1l}
\dot{x}_0=Ax_0+Bu_0,
\end{equation}
where $x_0\in\mathbf{R}^n$ is the state and
$u_0\in\mathbf{R}^{p}$ is the control input
of the leader. In some applications,
the leader might need its own control action
to achieve certain objectives,
e.g., to reach a desirable consensus value.
In this section, we consider the
general case where $u_0$ is possibly nonzero and time varying and
not accessible to any follower,
which is much harder to solve than the case with
$u_0=0$.

Before moving forward,
the following mild assumption is needed.

{\bf Assumption 4}~The leader's control input $u_0$ is bounded,
i.e., there exists a positive scalar $\gamma$ such that
$\|u_0\|\leq\gamma$.

It is assumed that the leader receives no information
from any follower and the state of the
leader is available to only a subset of the followers. The
communication graph among the $N+1$ agents is represented by a
directed graph $\widehat{\mathcal {G}}$, which satisfies the following
assumption.

{\bf Assumption 5}~$\widehat{\mathcal
{G}}$ contains a directed spanning tree with the leader as the root
and the subgraph associated with
the $N$ followers is undirected.

Denote by $\widehat{\mathcal {L}}$ the Laplacian matrix associated with
$\widehat{\mathcal {G}}$. Because the leader has no neighbors, $\widehat{\mathcal {L}}$
can be partitioned as
$\widehat{\mathcal {L}}=\begin{bmatrix} 0 & 0_{1\times N} \\
\mathcal {L}_2 & \mathcal {L}_1\end{bmatrix}$,
where $\mathcal {L}_2\in\mathbf{R}^{N\times 1}$ and $\mathcal
{L}_1\in\mathbf{R}^{N\times N}$.
By Lemma 1 and Assumption 5, it is clear that
$\mathcal {L}_1>0$.

The objective of this paper is to solve the leader-follower consensus
problem for the agents in \dref{1c} and \dref{1l}, i.e., to design
distributed protocols under which the states of the $N$ followers
converge to the state of the leader.

\subsection{Distributed Static Consensus Protocol}

Based on the relative states of neighboring agents,
the following distributed
static controller is proposed for each follower:
\begin{equation}\label{ssdf}
\begin{aligned}
u_i &=c_1K\sum_{j=0}^Na_{ij}(x_i-x_j)+[c_2+\rho_i(x_i,t)]\\
&\quad\times\tilde{g} (
K\sum_{j=0}^Na_{ij}(x_i-x_j) ),~i=1,\cdots,N,
\end{aligned}
\end{equation}
where $c_1>0$ and $c_2>0\in\mathbf{R}$ are constant coupling gains,
$a_{ij}$
is the $(i,j)$-th entry of the adjacency matrix
associated with $\widehat{\mathcal {G}}$,
the nonlinear function $\tilde{g}(\cdot)$ is
defined as follows:
for $w\in\mathbf{R}^n$,
\begin{equation}\label{satub}
\tilde{g}(w)=\begin{cases}\frac{w}{\|w\|}  & \text{if}~[\gamma +\rho_i(x_i,t)]\|w\|>\kappa\\
\frac{w}{\kappa} & \text{if}~[\gamma +\rho_i(x_i,t)]\|w\|\leq\kappa
\end{cases},
\end{equation}
with $\kappa$ being a small positive value,
and the rest of the variables
are defined as in \dref{ssd}.

Let $x=[x_1^T,\cdots,x_N^T]^T$,
$\zeta_i=x_i-x_0$, $i=1,\cdots,N$,
and $\zeta=[\zeta_1^T,\cdots,\zeta_N^T]^T$. Using \dref{ssdf} for
\dref{1c} and \dref{1l},
we can obtain the closed-loop network dynamics as
\begin{equation}\label{netss1f}
\begin{aligned}
\dot{\zeta}
&= (I_N\otimes A+c_1\mathcal {L}_1\otimes BK)\zeta+(I_N\otimes B)F(x,t)\\
&\quad+[(c_2 I_N+\rho(x,t))\otimes B]\widetilde{G}(\zeta)-({\bf
1}\otimes B)u_0,
\end{aligned}
\end{equation}
where $F(x,t)$ and $\rho(x,t)$ are defined as in \dref{netss1} and
\begin{equation}\label{sa3bf}
\widetilde{G}(\xi)\triangleq\begin{bmatrix}\tilde{g}(K\sum_{j=1}^N \mathcal {L}_{1j}\xi_j)\\
\vdots\\\tilde{g}(K\sum_{j=1}^N \mathcal {L}_{Nj}\xi_j)\end{bmatrix},
\end{equation}
with $\mathcal {L}_{ij}$
being the $(i,j)$-th entry of $\widehat{\mathcal {L}}$ associated with
$\widehat{\mathcal {G}}$.
Clearly, the leader-follower consensus problem is solved if
$\zeta$ of \dref{netss1} converges to zero. Hereafter,
we refer to $\zeta$ as the leader-follower consensus error.

{\bf Theorem 3}~Suppose that Assumptions 2, 4, and 5 hold.
The parameters in the distributed protocol \dref{ssdf}
are designed as $c_1\geq
\frac{1}{\lambda_{\min}(\mathcal
{L}_1)}$, $c_2\geq\gamma$,
and $K=-B^TP^{-1}$, where $P>0$ is a solution to the LMI \dref{alg1}.
Then, the leader-follower consensus error $\zeta$ of \dref{netss1f}
is uniformly ultimately bounded and exponentially
converges to the residual set
\begin{equation}\label{df1}
\mathcal {D}_4\triangleq \{\zeta : \|\zeta\|^2\leq\frac{2\lambda_{\max}(P)N\kappa
}{\alpha\lambda_{\min}(\mathcal
{L}_1)}\},
\end{equation}
with a convergence rate faster than
${\rm{exp}(-\alpha t)}$, where $\alpha$ is defined as in \dref{alf}.

{\bf Proof}~
Consider the following Lyapunov function candidate
$$V_3=\frac{1}{2}\zeta^T(\mathcal {L}_1\otimes P^{-1})\zeta.$$
By Lemma 1 and Assumption 5, we know that $\mathcal {L}_1>0$,
implying that $V_3$ is
positive definite. The time derivative of $V_3$
along the trajectory of \dref{netss1f} is given by
\begin{equation}\label{lyas2f}
\begin{aligned}
\dot{V}_3
&=\zeta^T(\mathcal {L}_1\otimes P^{-1}A+c_1\mathcal {L}_1^2\otimes P^{-1}BK)\zeta
\\
&\quad+\zeta^T[(\mathcal {L}_1\otimes P^{-1}B)F(x,t)-(\mathcal {L}_1{\bf
1}\otimes P^{-1}B)u_0]\\&\quad+\zeta^T[\mathcal {L}_1(c_2 I_N+\rho(x,t))\otimes P^{-1}B]\widetilde{G}(\zeta).
\end{aligned}
\end{equation}

In virtue of
Assumption 4, we have
\begin{equation}\label{lyas3f}
\begin{aligned}
-\zeta^T(\mathcal {L}_1{\bf
1}\otimes P^{-1}B)u_0
&\leq \sum_{i=1}^N \|B^TP^{-1}\sum_{j=1}^N\mathcal {L}_{ij}\zeta_j\|\|u_0\|\\
&\leq
\gamma \sum_{i=1}^N \|B^TP^{-1}\sum_{j=1}^N\mathcal {L}_{ij}\zeta_j\|.
\end{aligned}
\end{equation}
Similarly as in the proof of Theorem 1, it is easy to see that
\begin{equation}\label{lyas5f}
\begin{aligned}
\zeta^T(\mathcal {L}_1\otimes P^{-1}B) F(x,t)
\leq
 \sum_{i=1}^N\rho_i(x_i,t) \|B^TP^{-1}\sum_{j=1}^N\mathcal {L}_{ij}\zeta_j\|.
\end{aligned}
\end{equation}

Next, consider the following three cases.

i) $[\gamma +\rho_i(x_i,t)]\|K\sum_{j=1}^N \mathcal {L}_{ij}\zeta_j\|>\kappa$,
$i=1,\cdots,N$.

In this case, it follows from \dref{satu} and \dref{sa3bf} that
\begin{equation}\label{lyas4f}
\begin{aligned}
&\zeta^T[\mathcal {L}_1(c_2 I_N+\rho(x,t))\otimes P^{-1}B]\widetilde{G}(\zeta)
\\&\qquad=-\sum_{i=1}^N [c_2+\rho_i(x_i,t)]\|B^TP^{-1}\sum_{j=1}^N\mathcal {L}_{ij}\zeta_j\|.
\end{aligned}
\end{equation}
Substituting \dref{lyas3f}, \dref{lyas4f}, and \dref{lyas5f} into \dref{lyas2f} gives
$$\begin{aligned}
\dot{V}_3 &\leq
\frac{1}{2}\zeta^T\mathcal {Z}\zeta
-(c_2-\gamma)\sum_{i=1}^N \|B^TP^{-1}\sum_{j=1}^N\mathcal {L}_{ij}\zeta_j\|\\
&\leq \frac{1}{2}\zeta^T\mathcal {Z}\zeta,
\end{aligned}$$
where
$\mathcal {Z}=\mathcal {L}_1\otimes (P^{-1}A+A^TP^{-1})-2c_1\mathcal {L}_1^2\otimes P^{-1}BB^TP^{-1}$.

ii) $[\gamma +\rho_i(x_i,t)]\|K\sum_{j=1}^N \mathcal {L}_{ij}\zeta_j\|\leq\kappa$, $i=1,\cdots,N$.

In this case, we have
\begin{equation}\label{lyasc8f}
\begin{aligned}
&\zeta^T[\mathcal {L}_1(c_2 I_N+\rho(x,t))\otimes P^{-1}B]\widetilde{G}(\zeta)
\\& =-\sum_{i=1}^N \frac{c_2 +\rho_i(x_i,t)}{\kappa}\|B^TP^{-1}\sum_{j=1}^N \mathcal {L}_{ij}\zeta_j\|^2
\leq 0.
\end{aligned}
\end{equation}
Then, it follows from \dref{lyas2f}, \dref{lyas3f}, \dref{lyasc8f}, and \dref{lyas5f} that
\begin{equation}\label{lyasc9f}
\begin{aligned}
\dot{V}_3 &\leq \frac{1}{2}\xi^T\mathcal {Z}\zeta
+\sum_{i=1}^N [\gamma +\rho_i(x_i,t)]\|B^TP^{-1}\sum_{j=1}^N \mathcal {L}_{ij}\zeta_j\|\\
&\leq \frac{1}{2}\zeta^T\mathcal {Z}\zeta+N\kappa.
\end{aligned}
\end{equation}

iii)
$[\gamma +\rho_i(x_i,t)]\|K\sum_{j=1}^N \mathcal {L}_{ij}\zeta_j\|>\kappa$, $i=1,\cdots,l$, and
$[\gamma +\rho_i(x_i,t)]\|K\sum_{j=1}^N \mathcal {L}_{ij}\zeta_j\|\leq\kappa$, $i=l+1,\cdots,N$,
where $2\leq l\leq N-1$.
In this case, we can get that
\begin{equation}\label{lyasc10f}
\begin{aligned}
&\zeta^T[\mathcal {L}_1(c_2 I_N+\rho(x,t))\otimes P^{-1}B]\widetilde{G}(\zeta)
\\&\qquad\leq -\sum_{i=1}^l [c_2+\rho_i(x_i,t)]\|B^TP^{-1}\sum_{j=1}^N\mathcal {L}_{ij}\zeta_j\|.
\end{aligned}
\end{equation}
Then, it follows from \dref{lyas2f}, \dref{lyas3f}, \dref{lyasc10f}, and \dref{lyas5f} that
$$\begin{aligned}
\dot{V}_3 &\leq \frac{1}{2}\zeta^T\mathcal {Z}\zeta
-(c_2-\gamma)\sum_{i=1}^l \|B^TP^{-1}\sum_{j=1}^N\mathcal {L}_{ij}\zeta_j\|\\
&\quad+\sum_{i=l+1}^{N} [\gamma +\rho_i(x_i,t)]\|B^TP^{-1}\sum_{j=1}^N \mathcal {L}_{ij}\zeta_j\|\\
 &\leq \frac{1}{2}\zeta^T\mathcal {Z}\zeta+(N-l)\kappa.
\end{aligned}$$

Therefore, by examining the above three cases, we know that
$\dot{V}_3$ satisfies \dref{lyasc9f} for all $\zeta\in\mathbf{R}^{Mn}$.
Note that \dref{lyasc9f} can be rewritten as
\begin{equation}\label{lyasc11f}
\begin{aligned}
\dot{V}_3 &\leq
-\alpha V_3+\frac{1}{2}\zeta^T(\mathcal {Z}+\alpha\mathcal {L}_1\otimes P^{-1})\zeta+N\kappa.
\end{aligned}
\end{equation}
Because $\alpha= \frac{-\lambda_{\max}(AP+PA^T-2 BB^T)}{\lambda_{\max}(P)}$,
in light of \dref{alg1}, we can obtain that
$$\begin{aligned}
&(\mathcal {L}_1^{-\frac{1}{2}}\otimes P)
(\mathcal {Z}+\alpha\mathcal {L}_1\otimes P^{-1})
(\mathcal {L}_1^{-\frac{1}{2}}\otimes P)\\
&\quad=I_N\otimes (AP+PA^T+\alpha P)-2c_1\mathcal {L}_1\otimes BB^T\\
&\quad\leq I_N\otimes [AP+PA^T+\alpha P-2BB^T]<0.
\end{aligned}$$
Then, it follows from \dref{lyasc11f} that
$\dot{V}_3 \leq
-\alpha V_3+N\kappa$. The rest of the proof can be completed by
following similar steps in the proof of Theorem 1,
which is omitted here for conciseness.
\hfill $\blacksquare$

{\bf Remark 6}~ Different from \dref{ssd},
the term $c_2\tilde{g}(K\sum_{j=0}^Na_{ij}(x_i-x_j))$ in the consensus protocol \dref{ssdf} is used to deal with
the effect of the leader's nonzero control input $u_0$.
Due to the nonzero $u_0$, modifications are accordingly made
onto $g(\cdot)$ in \dref{satu} to get
the nonlinear function $\tilde{g}$ in \dref{sa3bf}. Similarly as in Theorem 1,
the parameters $c_1$, $c_2$, and $K$ are independently designed in Theorem 4,
where $c_2$ is related to only the upper bound of $u_0$.

\subsection{Distributed Adaptive Consensus Protocol}

In the last subsection, the design of the distributed protocol
\dref{ssdf} relies on $\lambda_{\min}(\mathcal {L}_1)$
and the upper bound $\gamma$ of the leader's
control input $u_0$ and the upper bounds $\rho_i(x_i,t)$ of
the matching uncertainties $f_i(x_i,t)$. The objective of this
subsection is to design a fully distributed
consensus protocol without requiring the aforementioned
global information. To this end,
we propose the following distributed adaptive
protocol to each follower as
\begin{equation}\label{sscaf}
\begin{aligned}
u_i &=\hat{d}_iK\sum_{j=0}^Na_{ij}(x_i-x_j)+r(
K\sum_{j=0}^Na_{ij}(x_i-x_j)),\\
\dot{\hat{d}}_i &= \tau_i [-\varphi_i\hat{d}_i+(\sum_{j=0}^Na_{ij}(x_i-x_j)^T)\Gamma
(\sum_{j=0}^Na_{ij}(x_i-x_j))\\
&\quad+\|K\sum_{j=0}^Na_{ij}(x_i-x_j)\| ],\\
\dot{\hat{e}}_i &=\epsilon_i [-\psi_i\hat{e}_i+\|K\sum_{j=0}^Na_{ij}(x_i-x_j)\|\|x_i\|],
i=1,\cdots,N,
\end{aligned}
\end{equation}
where
$\hat{d}_i(t)$ and $\hat{e}_i(t)$ are the adaptive gains
associated with the $i$-th follower and the rest of the variables
are defined as in \dref{ssca}.

Let the leader-follower consensus error
$\zeta$ be defined as in \dref{netss1f} and $\widehat{D}={\rm{diag}}(\hat{d}_1,\cdots,\hat{d}_N)$.
Then, it follows from \dref{1c}, \dref{1l}, and \dref{sscaf}
that the closed-loop network dynamics can be obtained as
\begin{equation}\label{netsca1f}
\begin{aligned}
\dot{\zeta}
&= (I_N\otimes A+\widehat{D}\mathcal {L}_1\otimes BK)\zeta+(I_N\otimes B)F(x,t)
\\&\quad+(I_N\otimes B)R(\zeta)-({\bf
1}\otimes B)u_0,\\
\dot{\hat{d}}_i &= \tau_i [-\varphi_i\bar{d}_i+(\sum_{j=1}^N\mathcal {L}_{ij}\zeta_j^T)\Gamma
(\sum_{j=1}^N\mathcal {L}_{ij}\zeta_j)+\|K\sum_{j=1}^N\mathcal {L}_{ij}\zeta_j\| ],\\
\dot{\hat{e}}_i &=\epsilon_i [-\psi_i\hat{e}_i+\|K\sum_{j=1}^N\mathcal {L}_{ij}\zeta_j\|\|x_i\| ],
~i=1,\cdots,N,
\end{aligned}
\end{equation}
where $R(\cdot)$ remains the same as in \dref{sa4}
and the rest of the variables are defined as in \dref{netss1f}.

To present the following theorem, we use a Lyapunov function in the form of
$$V_4=\frac{1}{2}\zeta^T(\mathcal {L}_1\otimes P^{-1})\zeta+\sum_{i=1}^N\frac{\breve{d}_i^2}{2\tau_i}
+\sum_{i=1}^N\frac{\breve{e}_i^2}{2\epsilon_i},$$
where $\breve{e}_i=\hat{e}_i-e_i$,
$\breve{d}_i=\hat{d}_i-\hat{\beta}$, $i=1,\cdots,N$,
and $\hat{\beta}\geq \max_{i=1,\cdots,N} \{d_i+\gamma,\frac{1}{\lambda_{\min}(\mathcal {L}_1)}\}$.

{\bf Theorem 4}~Supposing that Assumptions 3, 4, and 5 hold,
the leader-follower consensus error $\zeta$
and the adaptive gains $\hat{d}_i$
and $\hat{e}_i$, $i=1,\cdots,N$, in \dref{netsca1f}
are uniformly ultimately bounded under
the distributed adaptive protocol \dref{sscaf}
with $K$ and $\Gamma$ designed as in Theorem 3.
Moreover, the following two assertions hold.
\begin{itemize}
\item[i)]
For any $\varphi_i$ and $\psi_i$,
$\xi$, $\breve{d}_i$,
and $\breve{e}_i$ exponentially converge to the residual set
\begin{equation}\label{d5}
\mathcal {D}_5\triangleq \{\zeta,\breve{d}_i,\breve{e}_i:V_4<
\frac{1}{2\delta}\sum_{i=1}^N(\hat{\beta}^2\varphi_i+e_i^2\psi_i)+\frac{N\kappa}{4\delta}\},
\end{equation}
with a convergence rate faster than
${\rm{exp}(-\delta t)}$, where $ \delta$ and $\alpha$ are defined
as in Theorem 2 and \dref{alf}, respectively.
\item[ii)]
If $\varrho<\alpha$, where $\varrho$ is defined in Theorem 2,
then in addition to i),
$\zeta$ exponentially converges to the residual set
\begin{equation}\label{d7k}
\begin{aligned}
\mathcal {D}_6\triangleq \{\zeta:\|\zeta\|^2&\leq
\frac{\lambda_{\max}(P)}{\lambda_{\min}(\mathcal {L}_1)(\alpha-\varrho)}[\sum_{i=1}^N(\hat{\beta}^2\varphi_i+e_i^2\psi_i)
\\&\quad+\frac{1}{2}N\kappa]\},
\end{aligned}
\end{equation}
with a convergence rate faster than
${\rm{exp}(-\varrho t)}$
\end{itemize}

{\bf Proof}~It can be completed by following similar
steps as in the proofs of Theorems 2 and 3.
\hfill $\blacksquare$

{\bf Remark 7}~
In related works \cite{li2011trackingTAC,li2012adaptiveauto},
distributed protocols are designed
to achieve leader-follower consensus
for linear multi-agent systems with a
leader of bounded control input.
Compared to \cite{li2011trackingTAC,li2012adaptiveauto}
where the agents have linear nominal dynamics,
the multi-agent system considered in this section
is subject to nonidentical matching uncertainties,
for which case it is quite challenging to show the boundedness
of the consensus errors and the adaptive gains.
In \cite{Das20102014,zhang2012adaptive},
the distributed tracking problem of multi-agent
systems with unknown nonlinear dynamics are discussed,
where the agents are restricted
to be first-order and special high-order systems.
In contrast, this paper considers
general high-order multi-agent systems with matching uncertainties.
Contrary to the protocols
in \cite{Das20102014,zhang2012adaptive} whose design
depends on global information
of the communication graph, the adaptive consensus
protocols \dref{ssca} and \dref{sscaf}
proposed in this paper are fully distributed,
which do not require any global information.
Besides, both the cases with and without a leader
are addressed in this paper.


\section{Robustness With Respect To Bounded Non-matching Disturbances}

In the preceding sections, the external disturbances
in \dref{1c1} are assumed to satisfy
the matching condition, i.e., Assumption 1.
In this section, we examine the case
where the agents are subject to external disturbances
which do not necessarily satisfy
the matching condition and investigate
whether the proposed protocols in the preceding sections still ensure
the boundedness of the consensus error.
For conciseness, we consider here only
the case of the leaderless consensus problem. The case of the
leader-follower consensus problem can also be similarly discussed.

Consider a network of $N$ agents whose communication graph is
represented by an undirected graph $\mathcal
{G}$. The dynamics of the $i$-th agent is described by
\begin{equation}\label{1d1}
\dot{x}_i =Ax_i+B[u_i+f_i(x_i,t)]+\omega_i,
\end{equation}
where $f_i(x_i,t)$ is the lumped matching uncertainty defined as in \dref{1c}
and $\omega_i\in \mathbf{R}^{n}$
is the bounded non-matching external disturbance, satisfying

{\bf Assumption 6} There exist positive constants $\upsilon_i$ such that
$\|\omega_i\|\leq\upsilon_i$, $i=1,\cdots,N$.

First, we will investigate whether the distributed static protocol \dref{ssd}
can ensure the ultimate boundedness of the consensus error $\xi$ for the agents in \dref{1d1}.
Using \dref{ssd} for \dref{1d1},
we can obtain the closed-loop network dynamics in terms of $\xi$ as
\begin{equation}\label{netss2d}
\begin{aligned}
\dot{\xi}
&= (I_N\otimes A+c\mathcal {L}\otimes BK)\xi+(M\otimes B)F(x,t)\\
&\quad+[M\rho(x,t)\otimes B]G(\xi)+(M\otimes I)\omega,
\end{aligned}
\end{equation}
where $\omega=[\omega_1^T,\cdots,\omega_N^T]^T$
and the rest of the variables are defined as in \dref{netss2}.

The result redesign the static protocol \dref{ssd}
to guarantee the boundedness of $\xi$ in \dref{netss2d}.

{\bf Theorem 5}~
Suppose that the communication graph $\mathcal
{G}$ is connected
and Assumptions 2 and 5 hold. The consensus error
$\xi$ in \dref{netss2d} is ultimately bounded
under the static protocol \dref{ssd} with
$c\geq
\frac{1}{\lambda_2}$ and $K=-B^TQ^{-1}$, where
$Q>0$ is a solution to the following LMI:
\begin{equation}\label{alg2q}
AQ+QA^T+\varepsilon Q-2BB^T<0,
\end{equation}
where $\varepsilon>1$. Moreover,
$\xi$ exponentially converges to the residual
set
\begin{equation}\label{dd}
\mathcal {D}_7\triangleq \{\xi:\|\xi\|^2\leq\frac{2\lambda_{\max}(Q)
}{(\varepsilon-1)\lambda_2}[\frac{\lambda_{\max}(\mathcal {L})}{2\lambda_{\min}(Q)}\sum_{i=1}^{N}\upsilon_i^2+N\kappa]\},
\end{equation}
with a convergence rate faster than ${\rm{exp}(-(\varepsilon-1) t)}$.

{\bf Proof}.
Consider the following Lyapunov function candidate:
$$V_5=\frac{1}{2}\xi^T(\mathcal {L}\otimes Q^{-1})\xi.$$
By following similar steps as in the proof of Theorem 1, we can get
that the time derivative of $V_5$
along the trajectory of \dref{netss2d} satisfies
\begin{equation}\label{lyasc9d}
\begin{aligned}
\dot{V}_5 \leq \frac{1}{2}\xi^T\mathcal {R}\xi+\xi^T(\mathcal {L}\otimes Q^{-1})\omega+N\kappa,
\end{aligned}
\end{equation}
where
$\mathcal {R}\triangleq\mathcal {L}\otimes (Q^{-1}A+A^TQ^{-1})-2c \mathcal {L}^2\otimes Q^{-1}BB^TQ^{-1}.$
Using the following fact:
\begin{equation}\label{csq}
\begin{aligned}
&\frac{1}{2}\xi^T(\mathcal {L}\otimes Q^{-1})\xi-\xi^T(\mathcal {L}\otimes Q^{-1})\omega+\frac{1}{2}\omega^T(\mathcal {L}\otimes Q^{-1})\omega\\
&\qquad=\frac{1}{2}(\xi-\omega)^T(\mathcal {L}\otimes Q^{-1})(\xi-\omega)\geq0,
\end{aligned}
\end{equation}
we can get from \dref{lyasc9d} that
\begin{equation}\label{lyasc10d}
\begin{aligned}
\dot{V}_5
&\leq \frac{1}{2}\xi^T(\mathcal {R}+\mathcal {L}\otimes Q^{-1})\xi+\frac{1}{2}\omega^T(\mathcal {L}\otimes Q^{-1})\omega+N\kappa\\
&\leq \frac{1}{2}\xi^T(\mathcal {R}+\mathcal {L}\otimes Q^{-1})\xi+\frac{\lambda_{\max}(\mathcal {L})}{2\lambda_{\min}(Q)}\sum_{i=1}^{N}\upsilon_i^2+N\kappa.
\end{aligned}
\end{equation}
Note that \dref{lyasc10d} can be rewritten into
\begin{equation}\label{lyasc11d}
\begin{aligned}
\dot{V}_5
&\leq -(\varepsilon-1) V_5+\frac{1}{2}\xi^T[\mathcal {R}+\varepsilon\mathcal {L}\otimes Q^{-1}]\xi
\\&\quad+\frac{\lambda_{\max}(\mathcal {L})}{2\lambda_{\min}(Q)}\sum_{i=1}^{N}\upsilon_i^2+N\kappa.
\end{aligned}
\end{equation}
Letting $\bar{\xi}$ be defined as in the proof of Theorem 1,
we have
\begin{equation}\label{lyasc12d}
\begin{aligned}
&\xi^T[\mathcal {R}+\varepsilon\mathcal {L}\otimes Q^{-1}]\xi\\
&\qquad=\sum_{i=2}^{N}\lambda_i\bar{\xi}_i^T[AQ+QA^T+\varepsilon Q-2c\lambda_i BB^T]\bar{\xi}_i\\
&\qquad\leq \sum_{i=2}^{N}\lambda_i\bar{\xi}_i^T[AQ+QA^T+\varepsilon Q-2BB^T]\bar{\xi}_i\leq0.
\end{aligned}
\end{equation}
Therefore, we get from \dref{lyasc11d} and \dref{lyasc12d} that
\begin{equation}\label{lyasc13d}
\begin{aligned}
\dot{V}_5
&\leq -(\varepsilon-1) V_5
+\frac{\lambda_{\max}(\mathcal {L})}{2\lambda_{\min}(Q)}\sum_{i=1}^{N}\upsilon_i^2+N\kappa,
\end{aligned}
\end{equation}
which, together with \dref{lyas01}, implies
that $\xi$ exponentially converges to the residual set
$\mathcal {D}_7$ in \dref{dd} with a convergence rate not less than
${\rm{exp}(-(\varepsilon-1) t)}$.
\hfill $\blacksquare$

{\bf Remark 8}~
As shown in Proposition 1 in \cite{li2011dynamic},
there exists a $Q>0$ satisfying \dref{alg2q} if and only
if $(A,B)$ is controllable. Thus, a sufficient
condition for the existence of \dref{ssd} satisfying Theorem 5
is that $(A,B)$ is controllable, which, compared to the existence condition
of \dref{ssd} satisfying Theorems 1-4,
is stronger. 
It is worth mentioning
that large $\varepsilon$ in \dref{alg2q} yields a faster
convergence rate of the consensus error $\xi$,
but meanwhile generally implies
a high-gain $K$ in the protocol \dref{ssd}.
In implementation, a tradeoff has
to be made when choosing $\varepsilon$.

Next, we will redesign the distributed adaptive protocol \dref{ssca}
to ensure the boundedness of the consensus error
for the agents in \dref{1d1}.
Using \dref{ssca} for \dref{1d1},
we can obtain the closed-loop dynamics of the network as
\begin{equation}\label{netss2da}
\begin{aligned}
\dot{\xi}
&= (I_N\otimes A+M\overline{D}\mathcal {L}\otimes BK)\xi+(M\otimes B)F(x,t)
\\&\quad+(M\otimes B)R(\xi)+(M\otimes I)\omega,\\
\dot{\bar{d}}_i &= \tau_i [-\varphi_i\bar{d}_i+(\sum_{j=1}^N\mathcal {L}_{ij}\xi_j^T)\Gamma
(\sum_{j=1}^N\mathcal {L}_{ij}\xi_j)+\|K\sum_{j=1}^N\mathcal {L}_{ij}\xi_j\| ],\\
\dot{\bar{e}}_i &=\epsilon_i [-\psi_i\bar{e}_i+\|K\sum_{j=1}^N\mathcal {L}_{ij}\xi_j\|\|x_i\| ],
~ i=1,\cdots,N,
\end{aligned}
\end{equation}
where the variables are defined as in \dref{netss1}
and \dref{netss2d}.

{\bf Theorem 6}~Suppose that $\mathcal{G}$ is connected and
Assumptions 3 and 5 hold. Then,
both the consensus error $\xi$
and the adaptive gains $\bar{d}_i$
and $\bar{e}_i$, $i=1,\cdots,N$, in \dref{netsca1}
are uniformly ultimately bounded
under the distributed
adaptive protocol \dref{ssca}
with $K=-B^TQ^{-1}$ and $\Gamma=Q^{-1}BB^TQ^{-1}$, where
$Q>0$ is a solution to the LMI \dref{alg2q}.
Moreover, we have
\begin{itemize}
\item[i)]
For any $\varphi_i$ and $\psi_i$,
$\xi$, $\tilde{d}_i$,
and $\tilde{e}_i$ exponentially converge to the residual set
\begin{equation}\label{d2d}
\begin{aligned}
\mathcal {D}_8\triangleq \{\xi,\tilde{d}_i,\tilde{e}_i:V_6&<
\frac{1}{2\sigma}\sum_{i=1}^N(\beta^2\varphi_i+e_i^2\psi_i)\\
&\quad+\frac{\lambda_{\max}(\mathcal {L})}{2\sigma\lambda_{\min}(Q)}\sum_{i=1}^{N}\upsilon_i^2+\frac{N\kappa}{4\sigma}\},
\end{aligned}
\end{equation}
with a convergence rate faster than ${\rm{exp}(-\sigma t)}$,
where $\sigma\triangleq \min_{i=1,\cdots,N}\{\varepsilon-1,\varphi_i \tau_i,\psi_i\epsilon_i\}$
and
\begin{equation}\label{lyar11d}
V_6=\frac{1}{2}\xi^T(\mathcal {L}\otimes Q^{-1})\xi+\sum_{i=1}^N\frac{\tilde{d}_i^2}{2\tau_i}
+\sum_{i=1}^N\frac{\tilde{e}_i^2}{2\epsilon_i},
\end{equation}
where the variables are defined as in \dref{lyaas1}.
\item[ii)]
If $\varphi_i$ and $\psi_i$ satisfy
$\varrho<\varepsilon-1$, where $\varrho$ is defined in Theorem 2,
then in addition to i), $\xi$ exponentially converges to the residual set
\begin{equation}\label{d2kd}
\begin{aligned}
\mathcal {D}_{9}\triangleq \{\xi:\|\xi\|^2&\leq
\frac{\lambda_{\max}(Q)}{\lambda_2(\varepsilon-1-\varrho)}[\sum_{i=1}^N(\beta^2\varphi_i+e_i^2\psi_i)
\\&\quad+\frac{\lambda_{\max}(\mathcal {L})}{\lambda_{\min}(Q)}\sum_{i=1}^{N}\upsilon_i^2+\frac{1}{2}N\kappa]\}.
\end{aligned}
\end{equation}
\end{itemize}
with a convergence rate faster than ${\rm{exp}(-\varrho t)}$.

{\bf Proof}~
Choose the Lyapunov function candidate as in \dref{lyar11d}.
By following similar steps as in the proof of Theorem 2, it is not difficult to get that
the time derivative of $V_6$
along the trajectory of \dref{netss2da} can be obtained as
\begin{equation}\label{lyasacd5}
\begin{aligned}
\dot{V}_6
&\leq \frac{1}{2}\xi^T\mathcal {W}\xi+\xi^T(\mathcal {L}\otimes Q^{-1})\omega+\frac{1}{4}N\kappa\\
&\quad-\frac{1}{2}\sum_{i=1}^N(\varphi_i\tilde{d}_i^2+\psi_i\tilde{e}_i^2)
+\frac{1}{2}\sum_{i=1}^N(\beta^2\varphi_i+e_i^2\psi_i),
\end{aligned}
\end{equation}
where $\mathcal {W}\triangleq\mathcal {L}\otimes (Q^{-1}A+A^TQ^{-1})-2\beta \mathcal {L}^2\otimes Q^{-1}BB^TQ^{-1}.$
Using the fact \dref{csq},
we can get from \dref{lyasacd5} that
\begin{equation}\label{lyasacd6}
\begin{aligned}
\dot{V}_6
&\leq \frac{1}{2}\xi^T(\mathcal {W}+\mathcal {L}\otimes Q^{-1})\xi+\frac{1}{2}\omega^T(\mathcal {L}\otimes Q^{-1})\omega
+\frac{1}{4}N\kappa\\
&\quad-\frac{1}{2}\sum_{i=1}^N(\varphi_i\tilde{d}_i^2+\psi_i\tilde{e}_i^2)
+\frac{1}{2}\sum_{i=1}^N(\beta^2\varphi_i+e_i^2\psi_i),
\end{aligned}
\end{equation}
Note that \dref{lyasacd6} can be rewritten into
\begin{equation}\label{lyasacd61}
\begin{aligned}
\dot{V}_6
&\leq -\sigma V_6+\frac{1}{2}\xi^T[\mathcal {W}+(\sigma +1)\mathcal {L}\otimes Q^{-1}]\xi
\\
&\quad+\frac{1}{2}\omega^T(\mathcal {L}\otimes Q^{-1})\omega-\frac{1}{2}\sum_{i=1}^N[(\varphi_i-\frac{\sigma}{\tau_i})\tilde{d}_i^2\\
&\quad+(\psi_i-\frac{\sigma}{\epsilon_i})\tilde{e}_i^2)]
+\frac{1}{2}\sum_{i=1}^N(\beta^2\varphi_i+e_i^2\psi_i)+\frac{1}{4}N\kappa\\
&\leq -\sigma V_6+\frac{1}{2}\xi^T[\mathcal {W}+(\sigma +1)\mathcal {L}\otimes Q^{-1}]\xi\\
&\quad+\frac{\lambda_{\max}(\mathcal {L})}{2\lambda_{\min}(Q)}\sum_{i=1}^{N}\upsilon_i^2
+\frac{1}{2}\sum_{i=1}^N(\beta^2\varphi_i+e_i^2\psi_i)+\frac{1}{4}N\kappa,
\end{aligned}
\end{equation}
where we have used the fact that
$\sigma\leq\min_{i=1,\cdots,N}\{\varphi_i \tau_i,\psi_i\epsilon_i\}$
to get that last inequality.
Since $\sigma\leq\varepsilon-1$,
similarly as in the proof of Theorem 5,
we can show that $\xi^T[\mathcal {W}+(\sigma +1)\mathcal {L}\otimes Q^{-1}]\xi\leq0.$
Then, it follows from \dref{lyasacd61} that
\begin{equation}\label{lyasacd7}
\begin{aligned}
\dot{V}_6
&\leq
 -\sigma V_6+\frac{\lambda_{\max}(\mathcal {L})}{2\lambda_{\min}(Q)}\sum_{i=1}^{N}\upsilon_i^2
+\frac{1}{2}\sum_{i=1}^N(\beta^2\varphi_i+e_i^2\psi_i)\\
&\quad+\frac{1}{4}N\kappa,
\end{aligned}
\end{equation}
which implies that
$V_6$ exponentially converges to the residual set
$\mathcal {D}_8$ in \dref{d2d} with a convergence
rate faster than ${\rm{exp}(-\sigma t)}$.
By following similar steps as in the last part of the proof of Theorem 2,
it is not difficult to show that
for the case where $\varphi_i$ and $\psi_i$ satisfy
$\varrho<\varepsilon-1$,
$\xi$ exponentially converges to the residual set $\mathcal {D}_{9}$.
The details are omitted here for conciseness. \hfill $\blacksquare$

{\bf Remark 9}~
Compared to the residual sets of the consensus error $\xi$ for the agents in Sections 3 and 4,
the residual sets of $\xi$ for the agents in \dref{1d1}
further depends on the largest eigenvalue of
$\mathcal {L}$ and the magnitudes of the non-matching disturbances.
Contrary to the residual sets of
$\xi$ for the agents in \dref{1c1} which can be accordingly adjusted by properly choosing the design
parameters of the consensus protocols, the residual sets of
$\xi$ for the agents in \dref{1d1}
contain a constant term related to the magnitudes of the
non-matching disturbances.

\section{Simulation Examples}

In this section, two numerical examples are presented to illustrate
the theoretical results.

\begin{figure}[htbp]
\centering
\includegraphics[width=0.4\linewidth]{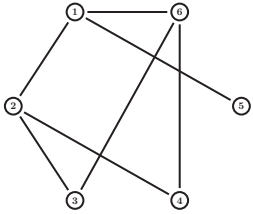}
\caption{The leaderless communication graph. }
\end{figure}

{\bf Example 1}~Consider a network of mass-spring systems
with a common mass $m$ but different unknown spring constants, described by
\begin{equation}\label{exa1}
m\ddot{y}_i+k_iy_i=u_i,~ i=1,\cdots,N,
\end{equation}
where $y_i$ are the displacements from certain reference positions and
$k_i$, $i=1,\cdots,N,$ are the bounded unknown spring constants.
Denote by $x_i=\left[\begin{smallmatrix} y_i& \dot{y}_i\end{smallmatrix}\right]^T$
the state of the $i$-th agent. Then, \dref{exa1} can be rewritten
as
\begin{equation}\label{exa13}
\dot{x}_i=Ax_i+B(u_i+k_iEx_i),~ i=1,\cdots,N,
\end{equation}
with $A=\left[\begin{smallmatrix} 0 & 1\\ 0 &0 \end{smallmatrix}\right]$, $B=\left[\begin{smallmatrix} 0 \\ \frac{1}{m}\end{smallmatrix}\right]$,
$E=\left[\begin{smallmatrix} -1 & 0\end{smallmatrix}\right].$ It is easy to see that
$k_iEx_i$, $i=1,\cdots,N$, satisfy Assumption 2, i.e., $\|k_iEx_i\|\leq k_i\|x_i\|$, $i=1,\cdots,N$.

Because the spring constants $k_i$ are unknown, we will use the adaptive protocol \dref{ssca}
to solve the consensus problem. Let $m=2.5kg$ and $k_i$ be randomly chosen.
Solving the LMI \dref{alg1} by using the Sedumi toolbox \cite{sturm1999using}
gives the feedback gain matrices of \dref{ssca} as
$K=-\left[\begin{smallmatrix} 0.6693  & 2.4595 \end{smallmatrix}\right]$ and
$\Gamma=\left[\begin{smallmatrix} 0.4480  &  1.6462\\
    1.6462  &  6.0489 \end{smallmatrix}\right].$
Assume that the communication topology is given in Fig. 1. In \dref{ssca},
select $\kappa=0.5$, $\varphi_i=\psi_i=0.05$,
and $\tau_i=\epsilon_i=10$, $i=1,\cdots,6$, in \dref{ssca}.
The state trajectories $x_i(t)$ of \dref{exa13} under
\dref{ssca} designed as above are depicted in Fig. 2,
which implies that consensus is indeed achieved.
The adaptive gains $\bar{d}_{i}$ and $\bar{e}_i$
in \dref{ssca} are shown in Fig. 3, from which it can be observed
that $\bar{d}_{i}$ and $\bar{e}_i$ tend to be quite small.

\begin{figure}[htbp]\centering
\includegraphics[width=0.45\linewidth,height=0.25\linewidth]{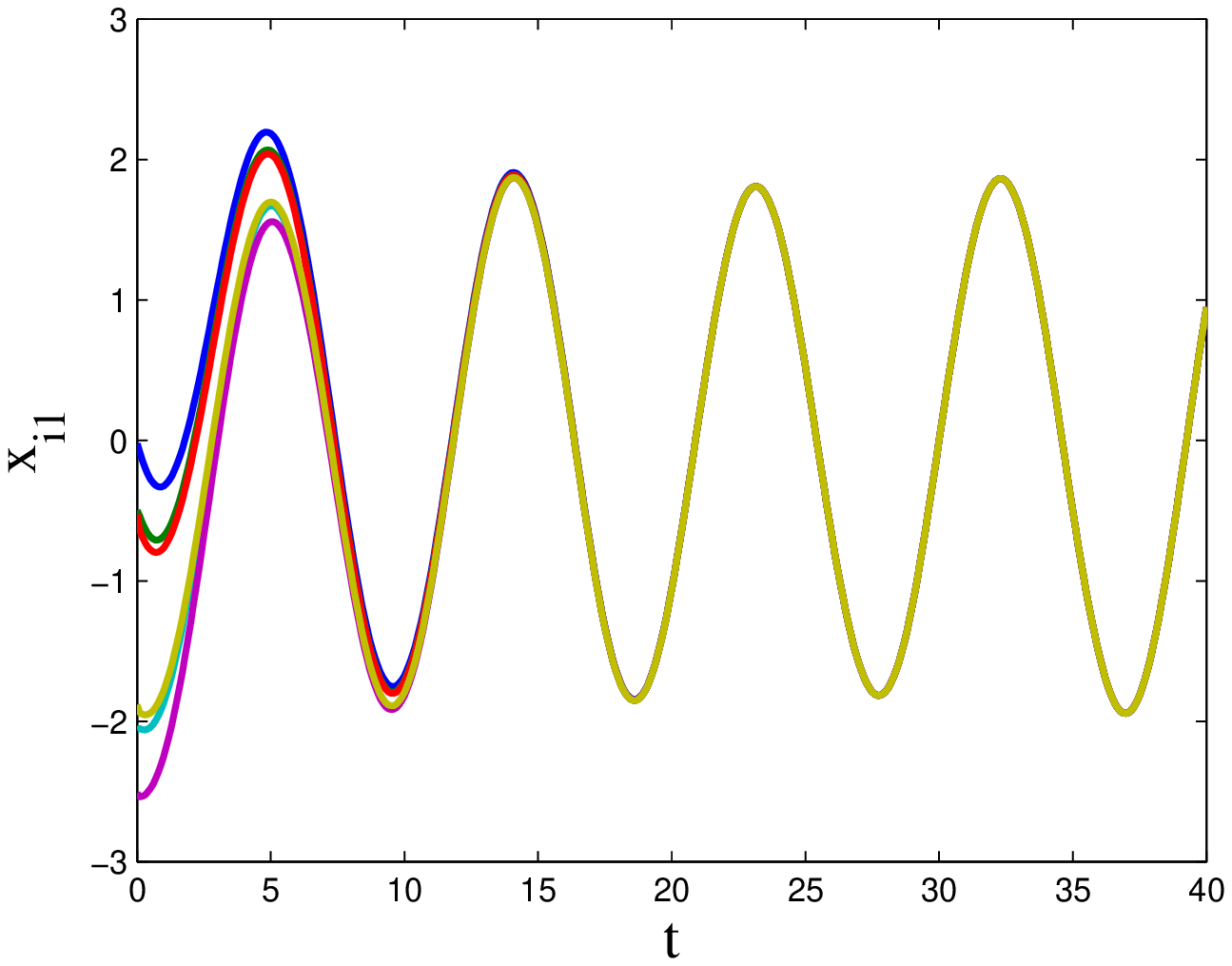}~
\includegraphics[width=0.45\linewidth,height=0.25\linewidth]{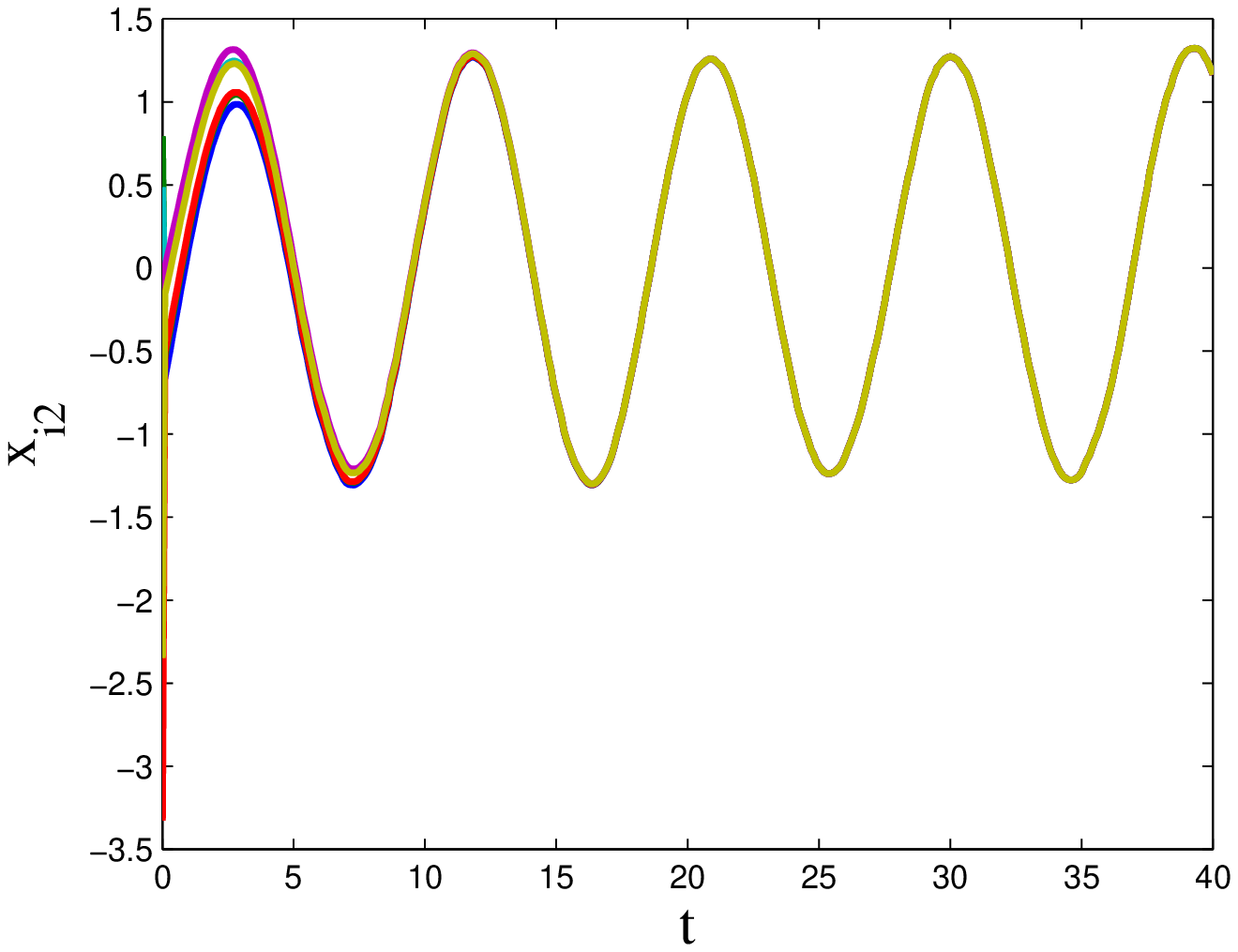}
\caption{The state trajectories of the mass-spring systems. }
\end{figure}

\begin{figure}[htbp]\centering
\centering
\includegraphics[width=0.45\linewidth,height=0.25\linewidth]{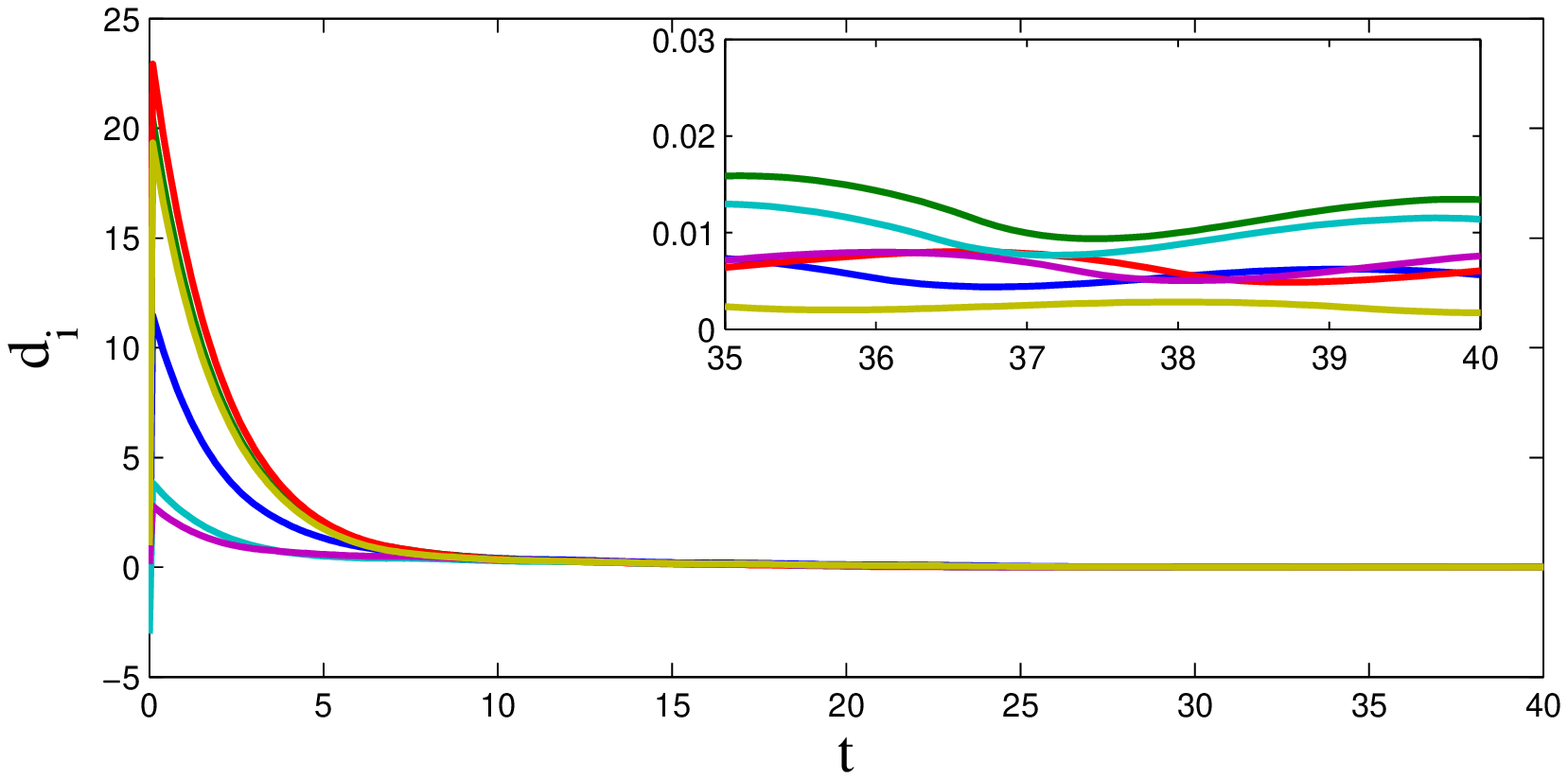}~
\includegraphics[height=0.25\linewidth,width=0.45\linewidth]{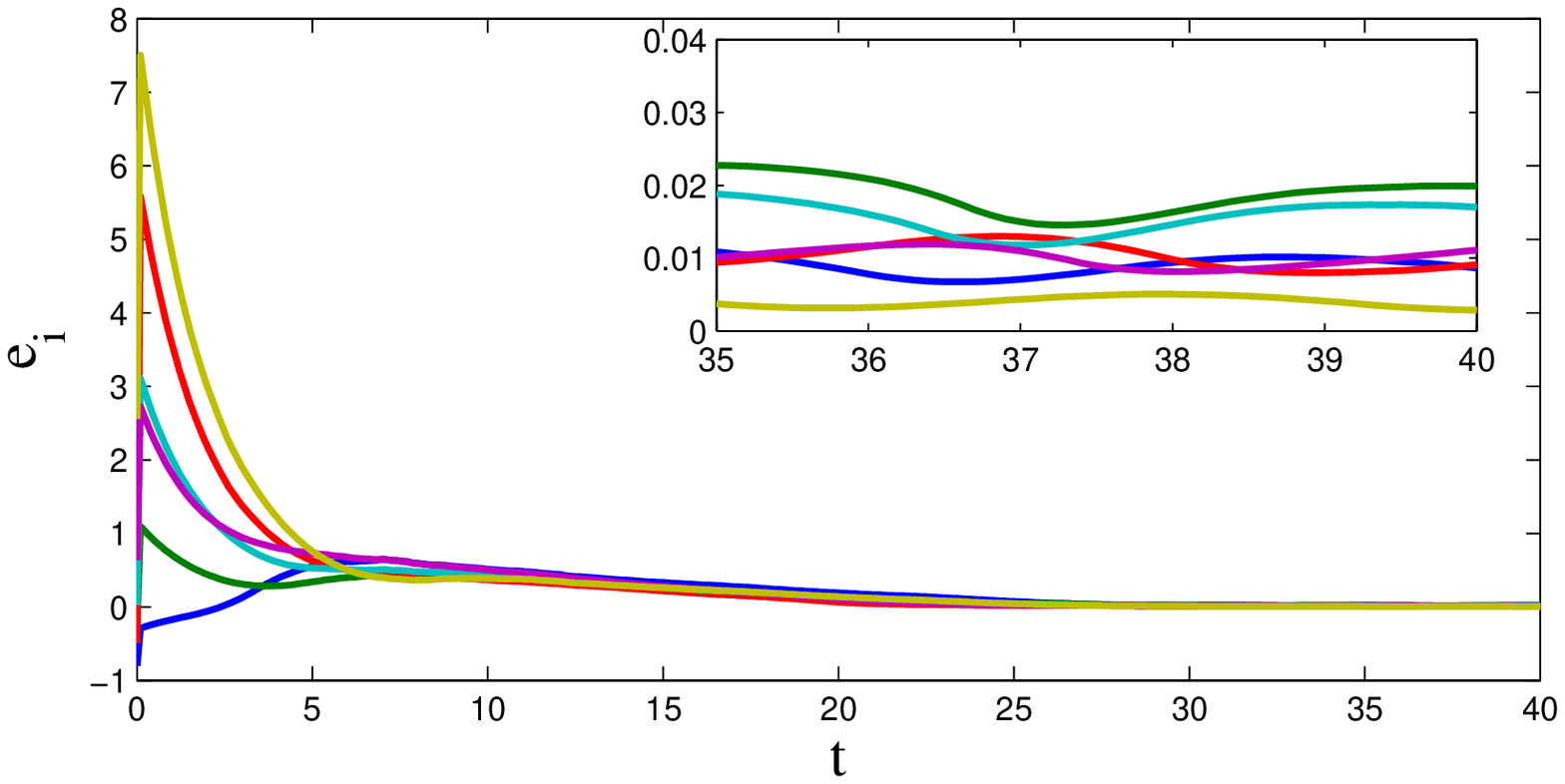}
\caption{The adaptive gains $\bar{d}_{i}$ and $\bar{e}_{i}$ in \dref{ssca}.}
\end{figure}

{\bf Example 2}~
Consider a group of Chua's circuits, whose dynamics in the dimensionless
form are given by \cite{madan1993chua}
\begin{equation}\label{chua}
\begin{aligned}
\dot{x}_{i1} &= a[-x_{i1}+x_{i2}-h(x_{i1})]+u_i,\\
\dot{x}_{i2} &=x_{i1}-x_{i2}+x_{i3},\\
\dot{x}_{i2} &= -bx_{i2},~ i=0,\cdots,N,
\end{aligned}
\end{equation}
where $a>0$, $b>0$, and $h(x_{i1})$ is a nonlinear function represented by
$h(x_{i1})=m_i^1 x_{i1}+\frac{1}{2}(m_i^2-m_i^1)(|x_{i1}+1|-|x_{i1}-1|),$
where $m_i^1<0$ and $m_i^2<0$.
The circuit indexed by 0 is the leader and the other circuits are the followers.
It is assumed that the Chua's circuits have nonidentical nonlinear components,
i.e., $m_i^1$ and $m_i^2$ are different for different Chua's circuits.
By letting $x_i=[x_{i1},x_{i2},x_{i3}]^T$, then \dref{chua} can be rewritten
in a compact form as
\begin{equation}\label{chua2}
\begin{aligned}
\dot{x}_{i} &= Ax_{i}+B[u_i+f_i(x_{i})],~ i=0,\cdots,N,
\end{aligned}
\end{equation}
where
$A=\left[\begin{smallmatrix}-m_0^1(a+1) & a & 0\\
1 & -1 &1 \\ 0 & -b & 0\end{smallmatrix}\right]$, $ B=\left[\begin{smallmatrix} 1 \\ 0 \\ 0\end{smallmatrix}\right]$,
$f_0(x_{0}) =\frac{a}{2}(m_0^1-m_0^2)(|x_{01}+1|-|x_{01}-1|)$,
$f_i(x_{i})=a(m_0^1-m_i^1) x_{i1}+\frac{a}{2}(m_i^1-m_i^2)(|x_{i1}+1|-|x_{i1}-1|)$, $i=1,\cdots,N$.
For simplicity, we let $u_0=0$ and
take $f_0(x_0)$ as the virtual control input of the leader,
which clearly satisfies $\|f_0(x_0)\|\leq\frac{a}{2}|m_0^1-m_0^2|$.
Let $a=9$, $b= 18$, $m_0^1=-\frac{3}{4}$, and $m_0^2=-\frac{4}{3}$.
In this case, the leader displays a
double-scroll chaotic attractor \cite{madan1993chua}.
The parameters $m_i^1$ and $m_i^2$, $i=1,\cdots,N$,
are randomly chosen within the interval $[-6,0)$.
It is easy to see that
$\|f_i(x_{i})\|\leq a|m_0^1-m_i^1|\|x_i\|+a|m_i^1-m_i^2|\\
=\frac{189}{4}\|x_i\|+54$, $i=1,\cdots,N.$
Note that $m_0^1$ is a parameter of the leader,
which might not be available to the followers.
Therefore, although $f_i(x_{i})$ satisfy
the above condition, the upper bound of $m_i^1-m_0^1$
might be not explicitly known for the followers.
Hence, we will use the adaptive protocol
\dref{sscaf} to solve the leader-follower consensus problem.

\begin{figure}[htbp]
\centering
\includegraphics[width=0.4\linewidth]{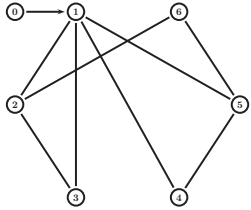}
\caption{The leader-follower communication graph. }
\end{figure}

\begin{figure}[htbp]\centering
\includegraphics[width=0.33\linewidth,height=0.25\linewidth]{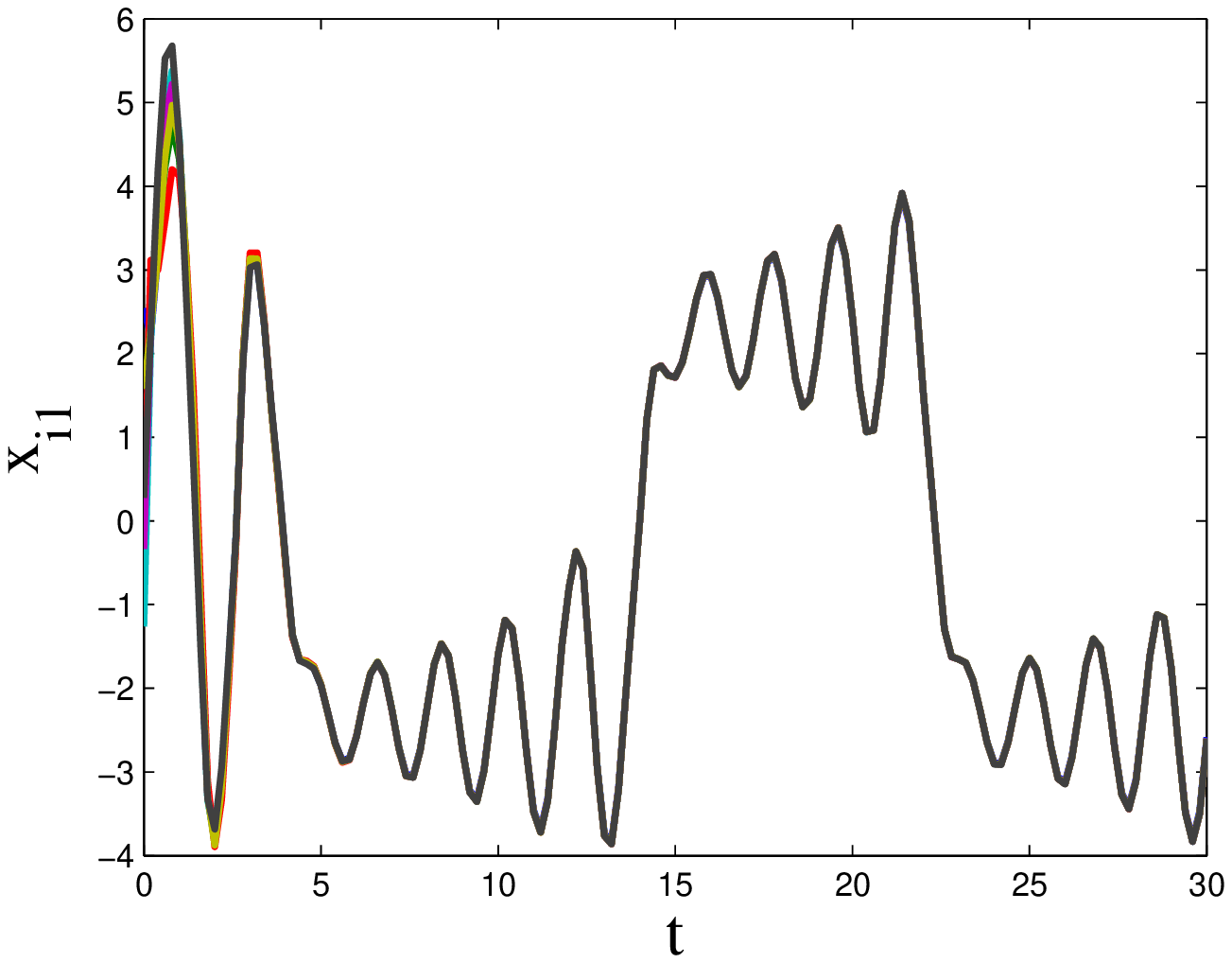}~
\includegraphics[width=0.33\linewidth,height=0.25\linewidth]{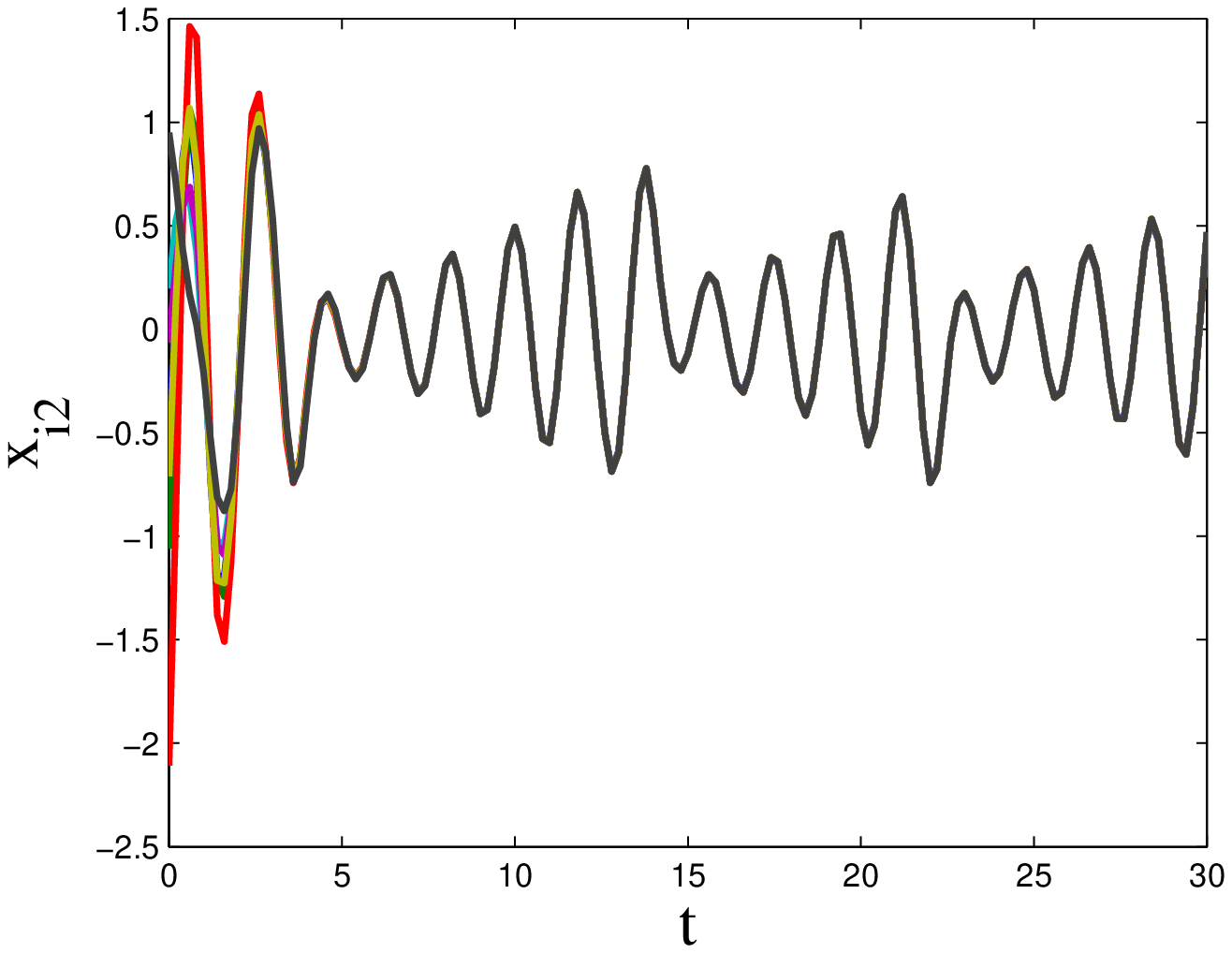}~
\includegraphics[width=0.33\linewidth,height=0.25\linewidth]{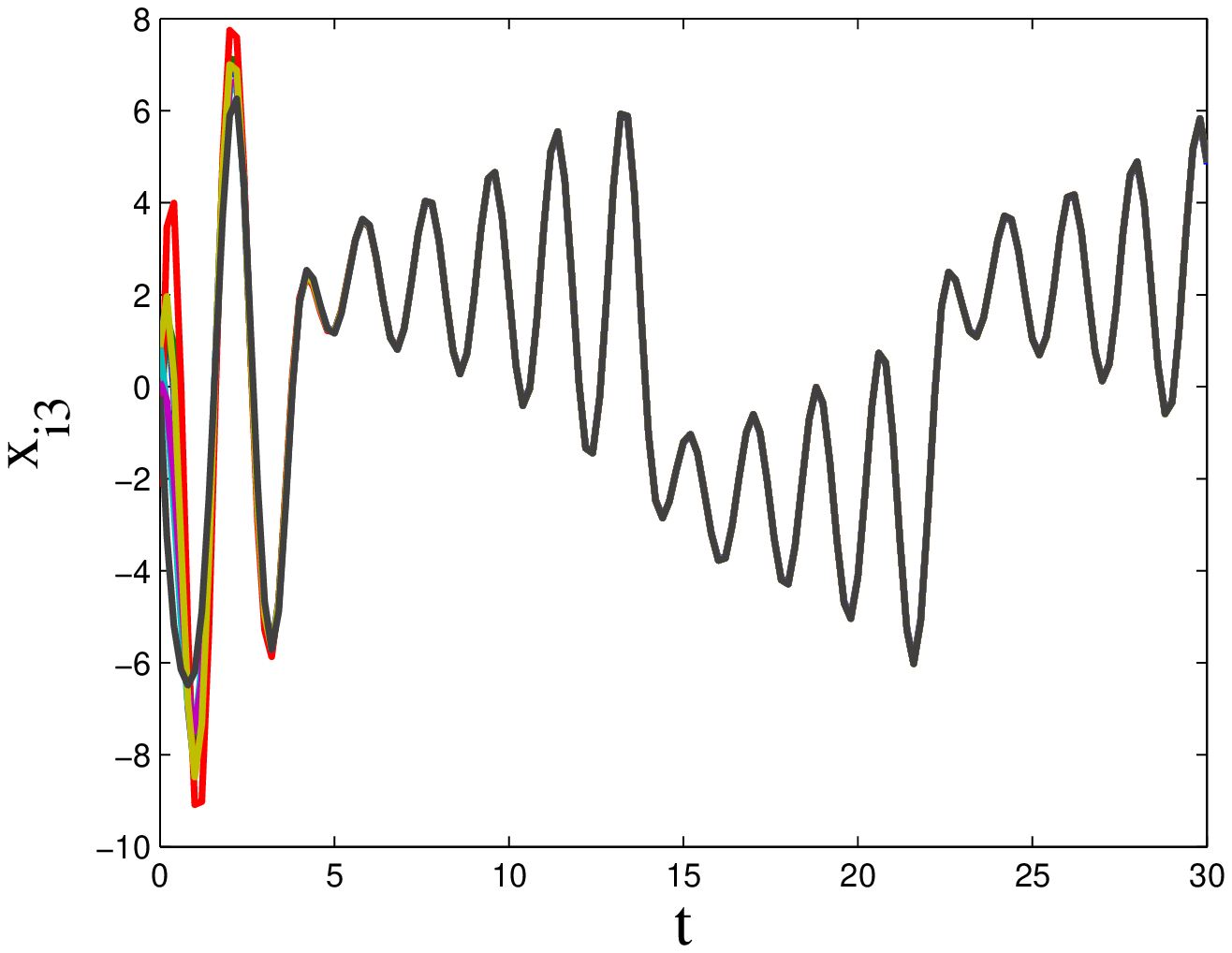}
\caption{The state trajectories of the Chua's circuits. }
\end{figure}

The communication
graph is given as in Fig. 4, where the node indexed by
0 is the leader. Solving the LMI \dref{alg1} gives the
feedback gain matrices of \dref{sscaf} as
$K=-\left[\begin{smallmatrix} 16.9070 & 16.5791  & 1.8297\end{smallmatrix}\right]$ and
$\Gamma=\left[\begin{smallmatrix} 285.8453 & 280.3016 &  30.9344\\
  280.3016 & 274.8654  & 30.3344\\
   30.9344 &  30.3344  &  3.3477\end{smallmatrix}\right].$
To illustrate Theorem 4,
select $\kappa=0.5$, $\varphi_i=\psi_i=0.05$,
and $\tau_i=\epsilon_i=5$, $i=1,\cdots,6$, in \dref{sscaf}.
The state trajectories $x_i(t)$ of the circuits under
\dref{sscaf} designed as above are depicted in Fig. 5,
implying that leader-follower consensus
is indeed achieved.
The adaptive gains
$\hat{d}_{i}$ and $\hat{e}_i$
in \dref{sscaf} are shown in Fig. 6, which
are clearly bounded.

\begin{figure}[htbp]\centering
\centering
\includegraphics[width=0.45\linewidth,height=0.25\linewidth]{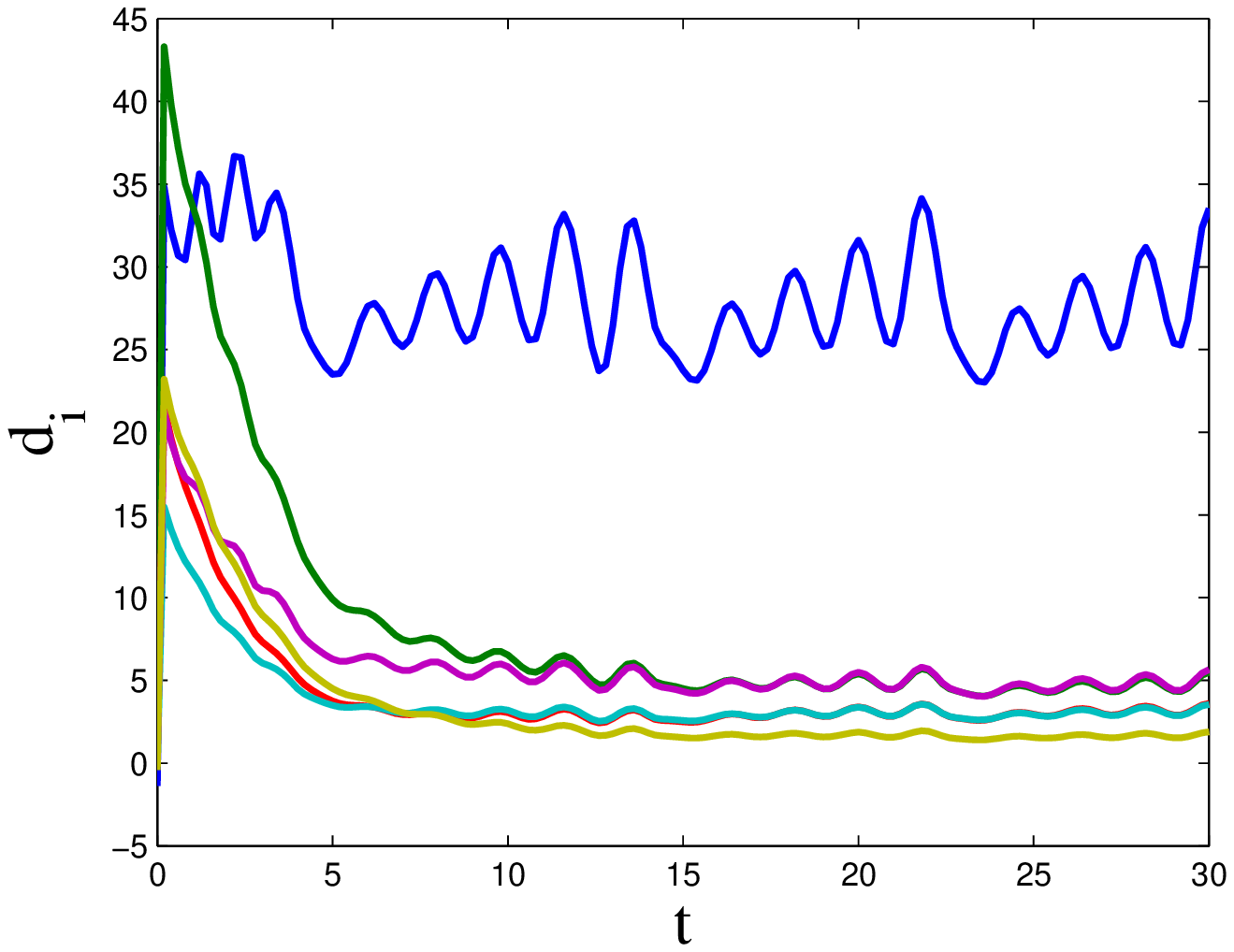}~
\includegraphics[height=0.25\linewidth,width=0.45\linewidth]{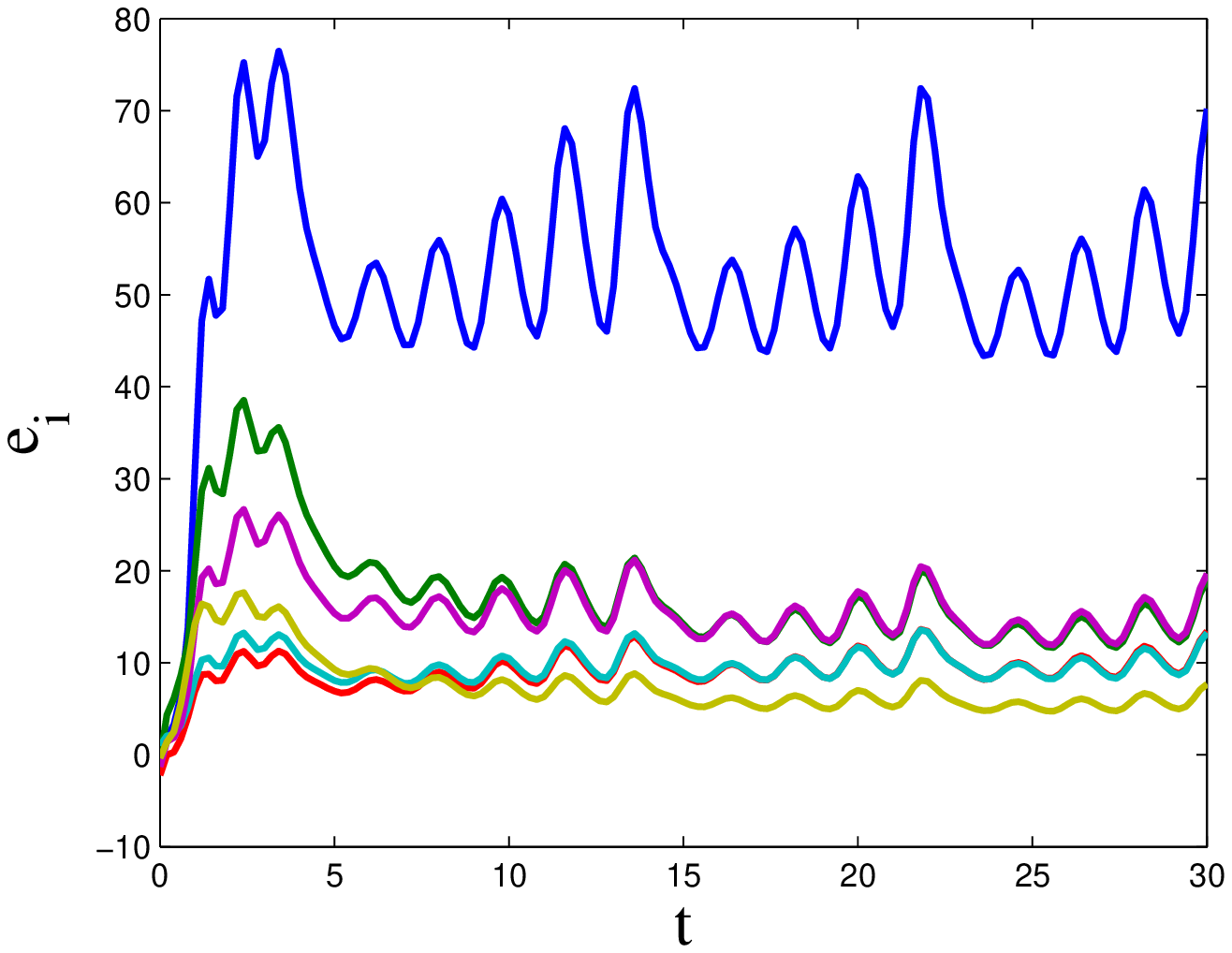}
\caption{The adaptive gains $\hat{d}_{i}$ and $\hat{e}_{i}$ in \dref{ssca}.}
\end{figure}

\section{Conclusion}

This paper has addressed the
robust consensus problem for
 multi-agent systems with heterogeneous matching uncertainties.
For both the cases with and without a
leader having a bounded unknown control input, several distributed continuous
static and adaptive consensus protocols have been designed, under which
the consensus error has been shown to be ultimately bounded and
exponentially converges to small adjustable residual sets. It should be noted that
the proposed adaptive consensus protocols can be implemented in a fully
distributed fashion without requiring any global information of
the communication graph or the upper bounds of the uncertainties and the leader's
control input. It has been also shown that the proposed protocols can be redesigned
so as to ensure the boundedness of the consensus error in the presence of
bounded external disturbances which do not necessarily satisfy the
matching condition. An interesting direction for future study
is to discuss the case with general directed
and switching communication graphs.


\end{document}